\documentclass[reqno,centertags]{amsart}
\usepackage{amsmath,amsthm,amscd,amssymb}
\usepackage{latexsym}
\newcommand{\bbC}{{\mathbb{C}}}
\newcommand{\bbD}{{\mathbb{D}}}

\newcommand{\bbH}{{\mathbb{H}}}

\newcommand{\bbR}{{\mathbb{R}}}
\newcommand{\bbT}{{\mathbb{T}}}
\newcommand{\bbU}{{\mathbb{U}}}
\newcommand{\bbZ}{{\mathbb{Z}}}

\newcommand{\calC}{{\mathcal{C}}}
\newcommand{\calE}{{\mathcal{E}}}

\newcommand{\calH}{{\mathcal H}}
\newcommand{\calL}{{\mathcal L}}
\newcommand{\calM}{{\mathcal M}}
\newcommand{\calS}{{\mathcal S}}

\newcommand{\dott}{\,\cdot\,}

\newcommand{\lb}{\label}
\newcommand{\f}{\frac}

\newcommand{\ol}{\overline}
\newcommand{\ti}{\tilde  }

\newcommand{\tr}{\text{\rm{Tr}}}
\newcommand{\dist}{\text{\rm{dist}}}

\newcommand{\Sz}{\text{\rm{Sz}}}

\newcommand{\spann}{\text{\rm{span}}}

\newcommand{\ran}{\text{\rm{ran}}}

\newcommand{\ac}{\text{\rm{ac}}}

\newcommand{\s}{\text{\rm{s}}}

\newcommand{\supp}{\text{\rm{supp}}}

\newcommand{\bi}{\bibitem}

\newcommand{\beq}{\begin{equation}}
\newcommand{\eeq}{\end{equation}}
\newcommand{\ba}{\begin{align}}
\newcommand{\ea}{\end{align}}
\newcommand{\veps}{\varepsilon}
\newcounter{smalllist}
\newenvironment{SL}{\begin{list}{{\rm\roman{smalllist})}}{%
\setlength{\topsep}{0mm}\setlength{\parsep}{0mm}\setlength{\itemsep}{0mm}%
\setlength{\labelwidth}{2em}\setlength{\leftmargin}{2em}\usecounter{smalllist}%
}}{\end{list}}

\DeclareMathOperator{\Real}{Re}

\DeclareMathOperator*{\wlim}{w-lim}
\DeclareMathOperator*{\esssup}{ess\,sup}

\allowdisplaybreaks
\numberwithin{equation}{section}
\newtheorem{theorem}{Theorem}[section]

\newtheorem*{p2.1}{Proposition 2.1}
\newtheorem{proposition}[theorem]{Proposition}
\newtheorem{lemma}[theorem]{Lemma}
\newtheorem{corollary}[theorem]{Corollary}
\theoremstyle{definition}

\theoremstyle{remark}

\newcommand{\abs}[1]{\lvert#1\rvert}

\begin{document}
\title{OPUC on One Foot}
\author[B. Simon]{Barry Simon*}

\thanks{$^*$ Mathematics 253-37, California Institute of Technology, Pasadena, CA 91125.
E-mail: bsimon@caltech.edu. Supported in part by NSF grant DMS-0140592}

\date{February 2, 2005} 

\begin{abstract} We present an expository introduction to 
orthogonal polynomials on the unit circle. 
\end{abstract}

\maketitle

%%%%%%%%%%%%%%%%%%%%%%%%%%%%%%%
\section{Introduction} \lb{s1} 
%%%%%%%%%%%%%%%%%%%%%%%%%%%%%%%

Orthogonal polynomials are the Rodney Dangerfield \cite{W} of analysis. Because of the  
impact of Stieltjes' great 1895 paper on F.~Riesz, Nevanlinna, and Hilbert's 
school, the moment problem and the closely related subject of orthogonal 
polynomials on the real line (OPRL) were central in the revolution in 
analysis from 1900--1920 and provided critical precursors to the Hahn-Banach 
theorem, the Riesz-Markov theorem, the spectral theorem, and the theory of 
selfadjoint extensions. But in recent years, too often the subject is dismissed 
as ``classical" and not worthy of further study.  

With developments in random matrix theory and combinatorics 
(e.g., \cite{BDK99,BDK00,BDMMZ,BR,BS,Joh00,MehBk,PrS00}), it is clear that orthogonal 
polynomials still have a lot to contribute. From one point of view, what makes 
them relevant is that they are the simplest of inverse spectral problems --- 
indeed, Gel'fand-Levitan \cite{GLev} explicitly note that their approach to 
inverse theory for Schr\"odinger operators is motivated by OPRL. Recently, 
OPUC ideas have provided a matrix realization of Lax pairs for the (defocusing) 
AKNS equation \cite{Nen}.

What is true for OPRL is even more true for orthogonal polynomials on the unit 
circle (OPUC). While the closely related area of positive harmonic functions 
on $\bbD=\{z\in\bbC\mid \abs{z}<1\}$ drew the attention of Carath\'eodory, 
Fej\'er, Herglotz, F.~Riesz, Schur, and Toeplitz in the 1910's, the subject was 
only invented by Szeg\H{o} in about 1920, especially in his deep 1920--1921 
paper \cite{Sz20-21}. So OPUC never had its era of centrality but has had a 
steady but small following over the years. Traditionally, the book references 
for the subject were Szeg\H{o}'s book \cite{Szb}, which has only one full 
and several partial chapters on OPUC, Geronimus' book \cite{GBk} and 
review \cite{GBk1}, and a chapter in Freud \cite{FrB}, which are very dated. 
With a major development published only in 2003 (the CMV matrix of 
Section~\ref{s5} below), it is hard not to be dated. Motivated by this 
dearth of review literature and by the opportunity to use Schr\"odinger 
operator techniques in a new setting, I published two volumes 
\cite{OPUC1,OPUC2} on the subject. Many friends asked if there wasn't 
some way to learn about the subject in less than 1100 pages, and this 
expository note is the result. 

Throughout, we use $\bbD$ for the unit disk in $\bbC$, and $\partial\bbD$ 
for the unit circle. Our inner products, $\langle f,g\rangle$, are linear in 
$g$ and antilinear in $f$.  Significant missing material involve some 
explicit examples --- these are discussed in Section~1.6 of \cite{OPUC1}: 
my favorite are the Rogers-Szeg\H{o} polynomials (Example~1.6.5). 
This article undergoes a kind of phase transition in the middle of 
Section~\ref{s5} in that before there, most results have proofs or at 
least sketches given, and afterwards, there aren't many proofs. This is 
because the earlier material is more central and also because the later 
proofs are lengthier. 

To put OPUC in context, recall some basics of OPRL. Since the fascinating 
issues of indeterminate moment problems (see \cite{AkhB,S270}) are 
irrelevant to OPUC, we will assume all measures have compact support: 
\begin{SL} 
\item[(1)] If $\mu$ is a probability measure on $\bbC$ (i.e., positive with 
$\mu(\bbC)=1$) with compact support and $X_n(z)$ are the monic orthogonal 
polynomials (i.e., $X_n(z)=z^n +$ lower order, $X_n\perp z^\ell$, $\ell 
=0,\dots, n-1$),  
\begin{align} 
zX_n(z) &= X_{n+1}(z) + \sum_{j=0}^n a_j^{(n)} X_j(z) \lb{1.1} \\
a_j^{(n)} &= \f{\langle X_j, zX_n\rangle}{\|X_j\|^2} \lb{1.2} 
\end{align} 
What makes OPRL special is that multiplication by $x$ is selfadjoint, so if 
we use $P_n$ in place of $X_n$ for OPRL and $\rho$ for $\mu$, 
\[
\langle P_j, xP_n\rangle = \langle xP_j, P_n\rangle =0 \qquad j=0, \dots, n-2 
\]
and thus \eqref{1.2} becomes 
\begin{equation} \lb{1.3} 
xP_n (x)=P_{n+1}(x) + b_{n+1} P_n(x) + a_n^2 P_{n-1}(x) 
\end{equation} 
for {\it Jacobi parameters}, $a_n, b_n$; $n=1,2,\dots$. If $p_n = P_n/\|P_n\|$ 
are the orthonormal OPRL, the matrix elements of multiplication by $x$ in 
$p_n$ basis have the form: 
\begin{equation} \lb{1.4} 
J=\begin{pmatrix} 
b_1 & a_1 & 0 & 0 & \dots \\
a_1 & b_2 & a_2 & 0 & \dots \\
0 & a_2 & b_3 & a_3 & \dots \\
\dots & \dots & \dots & \dots & \dots 
\end{pmatrix}
\end{equation} 

\item[(2)] There is a one-one correspondence between bounded $J$'s (i.e., $\sup_n 
\abs{a_n} + \abs{b_n}<\infty$) and $\rho$ on $\bbR$ with compact but infinite support. 
This is sometimes called Favard's theorem.

\item[(3)] If $A$ is a bounded selfadjoint operator on a separable Hilbert space, 
$\calH$, and $\varphi$ is a cyclic unit vector (i.e., $\{A^n\varphi\}_{n=0}^\infty$ 
span $\calH$), one can use the spectral theorem to find a measure $d\rho$ on 
$[-\|A\|,\|A\|]$ with $\int x^n\, d\rho =\langle \varphi, A^n\varphi\rangle$ and then 
the OPRL for this measure to find a semi-infinite Jacobi matrix unitarily equivalent 
to $A$ with  $\varphi$ mapped to $(1\, 0\, 0\, \dots)^t$. This realization is unique, 
that is, the $a_n$'s and $b_n$'s are intrinsic to the pair $(A,\varphi)$. It was Stone 
who emphasized this point of view that the study of Jacobi matrices was the same as 
the study of selfadjoint operators with a distinguished cyclic vector.  

\item[(4)] A key role is played by the Stieltjes transform of $\rho$, that is, 
the function, $m$, on $\bbC\backslash\supp(d\rho)$ given by 
\begin{equation} \lb{1.5} 
m(z) =\int \f{d\rho(x)}{x-z} 
\end{equation} 

\item[(5)] The Jacobi parameters can also be captured from $m(z)$ via a continued 
fraction expansion (of Stieltjes) at $\infty$: 
\begin{equation} \lb{1.6} 
m(z) = \cfrac{1}{-z+b_1 - \cfrac{a_1^2}{-z+b_2 - a_2^2 \dots} } 
\end{equation} 
\end{SL} 

\smallskip
We will not discuss applications of OPUC in detail but note its important applications 
to linear prediction and filtering theory. The basics are due to Wiener \cite{Wie}, 
Kolmogorov \cite{Kol}, Krein \cite{Krein1,Krein2}, and Levinson \cite{Lev}. The ideas 
have been especially developed by Kailath \cite{Kai74,Kai87,K6.1}. 

The title of this article is based on an incident reported in the Talmud \cite{Talmud} 
that someone asked the famous first-century rabbi Hillel to describe Judaism to him 
while he stood on one foot. Hillel's answer was: ``Do not do unto others that which 
is hateful to you. The rest is commentary. Go forth and study." This article is 
OPUC on one foot. \cite{OPUC1,OPUC2} are commentary. 

\smallskip 
It is a pleasure to thank M.~Aizenman for pushing me to write such an article. 
I'd like to thank S.~Denisov, F.~Gesztesy, L.~Golinskii, D.~Lubinsky, 
F.~Marcell\'an, P.~Nevai, and G.~Stolz for useful input. This paper was 
started while I was a visitor at the Courant Institute and completed during 
my stay as a Lady Davis Visiting Professor at Hebrew University, Jerusalem.  
I'd like to thank P.~Deift and C.~Newman for the hospitality of Courant and 
H.~Farkas and Y.~Last for the hospitality of the Mathematics Institute 
at Hebrew University.

%%%%%%%%%%%%%%%%%%%%%%%%%%%%%%%%%%%%%%%%%%
\section{The Szeg\H{o} Recursion} \lb{s2} 
%%%%%%%%%%%%%%%%%%%%%%%%%%%%%%%%%%%%%%%%%% 

OPUC is the study of probability measures on $\partial\bbD$, that is, positive 
measures, $\mu$, with 
\begin{equation} \lb{2.1} 
\mu(\partial\bbD)=1 
\end{equation} 
The {\it Carath\'eodory function} (after \cite{Ca07}) of $\mu$ is defined on 
$\bbD$ by 
\begin{equation} \lb{2.2} 
F(z) = \int \f{e^{i\theta}+z}{e^{i\theta}-z}\, d\mu(\theta)
\end{equation} 
This analog of \eqref{1.5} is an analytic function on $\bbD$ which obeys
\begin{equation} \lb{2.3} 
F(0)=1 \qquad z\in\bbD \Rightarrow \Real F(z) >0  
\end{equation} 
The {\it Schur function} (after \cite{Schur}) is then defined by 
\begin{equation} \lb{2.4} 
F(z) = \f{1+zf(z)}{1-zf(z)} 
\end{equation} 
and is an analytic function mapping $\bbD$ to $\ol{\bbD}$, that is, 
\begin{equation} \lb{2.4a} 
\sup_{z\in\bbD}\, \abs{f(z)} \leq 1 
\end{equation} 
($f(z)\equiv e^{i\theta_0}$ is included and produced by $\mu$, a point mass at 
$e^{i\theta_0}$). 

\eqref{2.2} sets up a one-one correspondence between probability measures $\mu$ and 
analytic functions obeying \eqref{2.3} --- this is essentially a form of the Herglotz 
representation (see \cite[pp.~247]{Rudin}) and can be realized via 
\begin{equation} \lb{2.5} 
d\mu =  \wlim_{r\uparrow 1}\, \Real F(re^{i\theta})\, \f{d\theta}{2\pi}
\end{equation} 
or by 
\begin{equation} \lb{2.6} 
F(z) = 1+2 \sum_{n=1}^\infty c_n z^n 
\end{equation} 
where $c_n$ are the {\it moments} of $\mu$ given by 
\begin{equation} \lb{2.7} 
c_n =\int e^{-in\theta} \, d\mu(\theta) 
\end{equation} 
\eqref{2.4} sets up a bijection between $f$'s obeying \eqref{2.4a} 
and $F$'s obeying \eqref{2.3}. 

We call a measure {\it trivial} if it is supported on a finite set and 
{\it nontrivial} otherwise. We will mainly be interested in nontrivial measures. 
$\mu$ is trivial if and only if its Schur function is a finite Blaschke product 
\begin{equation} \lb{2.8} 
f(z) = e^{i\theta_0} \prod_{j=1}^{n-1} \f{z-z_j}{1-\bar z_j z}
\end{equation} 
with $z_1, \dots, z_{n-1}\in\bbD$. Here $n$ is the number of points in the support 
of $d\mu$. Later (see the remark after Theorem~\ref{T7.1}) we will interpret 
\eqref{2.8} in terms of OPUC. 

If $\mu$ is a nontrivial probability measure on $\partial\bbD$, we define the 
monic orthogonal polynomials $\Phi_n (z;d\mu)$ (or $\Phi_n(z)$ if $d\mu$ 
is understood) by: 
\begin{equation} \lb{2.9} 
\Phi_n(z) =z^n + \text{ lower order} \quad \int e^{-ij\theta} 
\Phi_n (e^{i\theta}) \, d\mu(\theta) =0 \quad j=0,1,2,\dots, n-1
\end{equation} 
so in $L^2 (\partial\bbD, d\mu)$, $\langle \Phi_n, \Phi_m\rangle =0$ 
if $n\neq m$. The orthonormal polynomials $\varphi_n$ are defined by 
\begin{equation} \lb{2.10} 
\varphi_n(z) = \f{\Phi_n(z)}{\|\Phi_n\|} 
\end{equation} 
where $\|\cdot\|$ is the $L^2$-norm. $\{\varphi_n\}_{n=0}^\infty$ is an 
orthonormal set in $L^2$. It may not be a basis (e.g., $d\mu(\theta) = 
d\theta/2\pi$ where $\varphi_n(z) =z^n$ and $\ol{z^j}$, $j=1, \dots,$ are 
orthogonal to all $\varphi_n$). We will discuss this further below 
(see Theorem~\ref{T2.1A}).  

If $d\mu$ is trivial, say $\supp (d\mu) =\{z_j\}_{j=1}^k$, we 
can still define $\Phi_n, \varphi_n$ for $n=0,1,\dots, k-1$. We can even 
define $\Phi_k$ (but not $\varphi_k)$ as the unique monic polynomial of 
degree $k$ with $\|\Phi_k\|=0$, that is, 
\begin{equation} \lb{2.11} 
\Phi_k(z) =\prod_{j=1}^k (z-z_j) \qquad (\mu\text{ trivial})  
\end{equation} 

Clearly, \eqref{2.9} and the fact that the polynomials of degree at most $n$ 
have dimension $n+1$ implies 
\begin{equation} \lb{2.12} 
\deg(P) \leq n, \quad P\perp z^j, \quad j=0, \dots, n-1 \Rightarrow P=c\Phi_n
\end{equation} 

On $L^2 (\partial\bbD, d\mu)$, define the anti-unitary map, $^{*,n}$, by 
\begin{equation} \lb{2.13} 
f^{*,n} (e^{i\theta}) = e^{in\theta}\, \ol{f(e^{i\theta})} 
\end{equation} 
One mainly considers $^{*,n}$ on the set of polynomials of degree $n$ 
which is left invariant: 
\begin{equation} \lb{2.14} 
P(z) =\sum_{j=0}^n c_j z^j \Rightarrow P^{*,n}(z) = \sum_{j=0}^n \bar c_j  
z^{n-j} =z^n \, \ol{P(1/\bar z)} 
\end{equation}
Henceforth, following a standard, but unfortunate, convention, we drop the 
``$\, ,n$" and just use $P^*$, hoping the $n$ is implicit. Note that $1^* 
=z^n$, depending on $n$! 

Since $^*$ is anti-unitary, \eqref{2.12} implies 
\begin{equation} \lb{2.15} 
\deg(P)\leq n, \quad P\perp z^j, \quad j=1, \dots, n \Rightarrow P=c\Phi_n^* 
\end{equation} 
Since $\langle f,zg\rangle =\langle z^{-1}f, g\rangle$, it is easy to see that 
$\Phi_{n+1} -z\Phi_n \perp z^j$ for $j=1,2, \dots, n$. Since $\Phi$ is monic, 
this difference is of degree $n$, so \eqref{2.15} implies 
\begin{equation} \lb{2.16} 
\Phi_{n+1}(z) =z\Phi_n(z) -\bar\alpha_n \Phi_n^*(z) 
\end{equation} 
for complex numbers $\alpha_n$, called the {\it Verblunsky coefficients} 
(in the older literature, also called reflection, Schur, Szeg\H{o}, or 
Geronimus coefficients). \eqref{2.16} is called {\it Szeg\H{o} recursion} 
after its first occurrence in Szeg\H{o}'s book \cite{Szb}. In the engineering 
literature, it is called the Levinson algorithm after its rediscovery in 
linear prediction theory \cite{Lev}. The choice of minus and $\bar\alpha_n$ 
rather than $\alpha_n$ will be made clear by Geronimus' theorem (see 
Theorem~\ref{T3.1}). Since $\Phi_n$ is monic, \eqref{2.14} implies $\Phi_n^*(0) 
=1$, so \eqref{2.16} at $z=0$ implies 
\begin{equation} \lb{2.16a} 
\alpha_n =-\ol{\Phi_{n+1}(0)} 
\end{equation} 

\begin{theorem}\lb{T2.1} We have 
\begin{align} 
\|\Phi_{n+1}\|^2 &= (1-\abs{\alpha_n}^2) \|\Phi_n\|^2 \lb{2.17} \\
\|\Phi_n\| &= \prod_{j=0}^{n-1} (1-\abs{\alpha_j}^2)^{1/2} \lb{2.18}
\end{align}
For any nontrivial $\mu$, we have $\alpha_j (d\mu)\in\bbD$ for all $j$. 
If $\mu$ is trivial with $n$ points in its support, then $\alpha_0 (d\mu), 
\dots, \alpha_{n-2} (d\mu)\in\bbD$ and $\alpha_{n-1}(d\mu)\in \partial\bbD$. 
\end{theorem} 

\begin{proof} \eqref{2.16}, unitarity of multiplication by $z$, and $\Phi_n^* 
\perp\Phi_{n+1}$ imply 
\[
\|\Phi_n\|^2 = \|z\Phi_n\|^2 = \|\Phi_{n+1} + \bar\alpha_n \Phi_n^*\|^2 
= \|\Phi_{n+1}\|^2 + \abs{\alpha_n}^2 \|\Phi_n\|^2 
\]
which implies \eqref{2.17}. Induction and $\Phi_0 =1$ implies \eqref{2.18}. 
By \eqref{2.17}, $\abs{\alpha_j} <1$ in the nontrivial case and for $j=0, 
\dots, n-2$ in the trivial case. Since $\|\Phi_n\|=0\neq \|\Phi_{n-1}\|$ in 
the trivial case, \eqref{2.17} implies $\abs{\alpha_{n-1}}=1$. 
\end{proof} 

Since it arises often, we define 
\begin{equation} \lb{2.19} 
\rho_j =(1-\abs{\alpha_j}^2)^{1/2} \qquad \abs{\alpha_j}^2 + \rho_j^2 =1 
\end{equation} 

One can use \eqref{2.18} to relate completeness of $\{\varphi_n\}_{n=0}^\infty$ to 
the Verblunsky coefficients: 

\begin{theorem}\lb{T2.1A} For any nontrivial measure, the following are equivalent: 
\begin{SL} 
\item[{\rm{(a)}}] $\lim_{n\to\infty} \|\Phi_n\| =0$ 
\item[{\rm{(b)}}] $\sum_{j=0}^\infty \abs{\alpha_j}^2 =\infty$ 
\item[{\rm{(c)}}] $\{\varphi_n\}_{n=0}^\infty$ are a basis for $L^2 (\partial\bbD, 
d\mu)$ 
\end{SL} 
\end{theorem} 

{\it Remark.} We will see later that there is an additional equivalence via 
Szeg\H{o}'s theorem (see \eqref{8.8}). The equivalence of a Szeg\H{o} condition to 
completeness is due to Kolmogorov \cite{Kol} and Krein \cite{Krein1,Krein2}. 

\begin{proof}[Sketch] By \eqref{2.18}, (a) $\Leftrightarrow$ (b). If 
\begin{equation} \lb{2.18a} 
P_{[k,\ell] } = \text{projection in } L^2 (\partial\bbD, d\mu) 
\text{ onto } \spann\{z^m\}_{m=k}^\ell
\end{equation} 
we have that 
\begin{align}
\|\Phi_n\| &= \|(1-P_{[0,n-1]})z^n\|  \lb{2.18b} \\
&= \|(1-P_{[1,n]})1\| \lb{2.18c} \\
&= \|(1-P_{[0,n-1]})z^{-1}\| \lb{2.18d}
\end{align} 
where \eqref{2.18b} follows from the definition of $\Phi_n$, \eqref{2.18c} 
by applying $\,^{*,n}$ to $z^n$ and $P_{[0,n-1]}$, and \eqref{2.18d} by using 
the fact that multiplication by $z^{-1}$ is unitary. It follows that 
\begin{equation} \lb{2.19e} 
\|(1-P_{[0,\infty)})z^{-1}\| =\lim_{n\to\infty}\, \|\Phi_n\| 
\end{equation}
so (a) $\Leftrightarrow z^{-1}\in\spann\{\varphi_n\}_{n=0}^\infty$. If $z^{-1} 
\notin\spann\{\varphi_n\}_{n=0}^\infty$, clearly they are not complete. If 
$z^{-1}\in\spann\{\varphi_j\}_{j=0}^\infty$, an argument (see the proof of 
Theorem~1.5.7 in \cite{OPUC1}) taking powers of $z^{-1}$ shows $z^{-\ell} 
\in\spann \{\varphi_n\}_{n=0}^\infty$ for all $\ell$, so 
$\{\varphi_n\}_{n=0}^\infty$ are complete.  
\end{proof} 

Let $\bbD^{\infty, c}$ denote the set of complex sequences $\{\alpha_j\}_{j=0}^N$ 
where either $N=\infty$ and $\abs{\alpha_j}<1$ for all $j$, or else $N<\infty$ and 
$\alpha_0, \dots, \alpha_{N-1}\in\bbD$ while $\alpha_N\in\partial\bbD$. In the 
topology of componentwise convergence, $\bbD^{\infty,c}$ is compact (and is a 
compactification of $\bbD^\infty$). The map, $\calS$, from $\mu\mapsto 
\{\alpha_j(d\mu)\}_{j=0}^N$ is a well-defined map from $\calM_{+,1}(\partial\bbD)$, 
the probability measures on $\partial\bbD$, to $\bbD^{\infty,c}$. By \eqref{2.16}, 
the $\alpha$'s determine the $\Phi_n$'s. Since $\int \Phi_n(z)\, d\mu = \delta_{n0}$, 
the $\Phi_n$'s determine the moments inductively, and so $d\mu$, since 
$\{z^\ell\}_{\ell=-\infty}^\infty$ span a dense set of $C(\partial\bbD)$. 
Thus $\calS$ is one-one. Moreover, 

\begin{theorem}[Verblunsky's Theorem \cite{V35}] \lb{T2.2} $\calS$ is onto. 
\end{theorem} 

\cite{OPUC1} has four proofs of this theorem (Theorems~1.7.11, 3.1.3, 4.1.5, and 
4.2.8); see Section~\ref{s3} below. Given that $\calS$ is a bijection, it is easy to see 
that it is a homeomorphism if $\calM_{+,1}(\partial\bbD)$ is given the 
vague (i.e., $C(\partial\bbD)$-weak $^*$) topology.  

Applying $^*$ (actually, $^{*,n+1}$) to \eqref{2.16} yields 
\begin{equation} \lb{2.19a} 
\Phi_{n+1}^*(z) = \Phi_n^*(z) -\alpha_n z \Phi_n(z) 
\end{equation} 
Using \eqref{2.17} and \eqref{2.10}, we get the recursion relations for $\varphi_n$ 
written in matrix form
\begin{equation} \lb{2.19b} 
\binom{\varphi_{n+1}(z)}{\varphi_{n+1}^*(z)} = A(z,\alpha_n) 
\binom{\varphi_n(z)}{\varphi_n^*(z)} 
\end{equation}
where 
\begin{equation} \lb{2.20} 
A(z,\alpha) = \rho^{-1} 
\begin{pmatrix} z & -\bar\alpha \\ -z\alpha & 1 \end{pmatrix}
\end{equation}
Notice that $\det A=z$, so by inverting $A$, we get inverse recursion relations. 
We note the one for $\Phi_{n-1}$: 
\begin{equation} \lb{2.21} 
\Phi_{n-1}(z) = \f{\rho_{n-1}^{-2} [\Phi_n + \bar\alpha_{n-1}\Phi_n^*]}{z} 
\end{equation} 
Note that, by \eqref{2.16a}, $[\dots]$ vanishes at zero, so the right side of 
\eqref{2.21} is a polynomial of degree $n-1$. This implies: 

\begin{theorem}[Geronimus \cite{Ger46}]\lb{T2.3} Let $\mu,\nu$ be two probability 
measures on $\partial\bbD$ so that for some $N_0$, $\Phi_{N_0}(z;d\mu) = 
\Phi_{N_0}(z;d\nu)$. Then $\Phi_j (z;d\mu) =\Phi_j (z;d\nu)$ for $j=0,1,\dots, 
N_0 -1$, $\alpha_j (d\mu) =\alpha_j (d\nu)$ for $j=0,1,\dots, N_0-1$, and $\varphi_j 
(z;d\mu)= \varphi_j (z;d\nu)$ for $j=0,1,\dots, N_0$. 
\end{theorem} 

{\it Remark.} As noted in a footnote in Geronimus \cite{Ger46} and rediscovered 
by Wendroff \cite{Wen}, the result for OPRL requires equality for $P_{N_0}$ 
and $P_{N_0-1}$ and, in particular, it often happens that $P_{N_0}(x,d\gamma) = 
P_{N_0}(x,d\rho)$, but no other $P_j$'s are equal. 

\begin{proof} By \eqref{2.16a}, $\Phi_{N_0}$ at $0$ determines $\alpha_{N_0-1}$, 
and so $\rho_{N_0-1}$, and thus $\Phi_{N_0-1}$ by \eqref{2.21}. By induction, 
all $\alpha_j$, $j\leq N_0-1$, and $\Phi_j$, $j\leq N_0$, are equal and so, 
by \eqref{2.18}, $\|\Phi_j\|$, and so $\varphi_j$. 
\end{proof} 

As a final aspect of Szeg\H{o} recursion, we turn to the Christoffel-Darboux 
formula (proven by Szeg\H{o} \cite{Szb} for OPUC; Christoffel \cite{Chris} and 
Darboux \cite{Dar} had a similar formula for OPRL), which is an analog of an 
iterated Wronskian formula for ODE's. With $A$ given by \eqref{2.20}, one finds, 
by matrix multiplication, that 
\begin{equation} \lb{2.30a}
A(\zeta, \alpha_n)^* \begin{pmatrix} -1 & 0 \\ 0 & 1 \end{pmatrix} 
A(z,\alpha_n) = \begin{pmatrix} -z\bar\zeta & 0 \\ 0 & 1 \end{pmatrix} 
\end{equation}
so 
\begin{align*} 
\ol{\varphi_{n+1}^*(\zeta)}\, \varphi_{n+1}^*(z) &- 
\ol{\varphi_{n+1} (\zeta)}\, \varphi_{n+1}(z) \\
&= \ol{\varphi_n^*(\zeta)}\, \varphi_n^*(z) - z\bar\zeta \,\, \ol{\varphi_n(\zeta)}\, 
\varphi_n(z) \\ 
&= (1-z\bar\zeta) \, \ol{\varphi_n(\zeta)}\, \varphi_n(z) + 
[\, \ol{\varphi_n^*(\zeta)}\, \varphi_n^*(z) - \ol{\varphi_n^*(\zeta)}\, 
\varphi_n(z)] 
\end{align*} 
which, iterated to $n=0$ (where $[\dots]=0$), yields 

\begin{theorem}[Szeg\H{o} \cite{Szb}; CD Formula for OPUC]\lb{T2.4} 
\begin{equation} \lb{2.21x} 
(1-z\bar\zeta) \sum_{j=0}^n \, \ol{\varphi_n(\zeta)}\, \varphi_n(z) = 
\ol{\varphi_{n+1}^*(\zeta)}\, \varphi_{n+1}^*(z) - 
\ol{\varphi_{n+1}(\zeta)}\, \varphi_{n+1}(z)
\end{equation} 
\end{theorem}

If $z=\zeta$ and lie in $\bbD$, we have various positivity facts that 
imply {\rm{(}}the first since $\varphi_0(z)=1${\rm{)}}: 
\begin{alignat}{2}
\abs{\varphi_n^*(z_0)} &\geq (1-\abs{z_0}^2)^{1/2}&&\qquad\text{for } 
z_0\in\bbD  \lb{2.22} \\ 
\lim\, \abs{\varphi_{n+1}^*(z_0)} &= \infty\Leftrightarrow \sum_{j=0}^\infty 
\abs{\varphi_j (z_0)}^2=\infty &&\qquad \text{for } z_0\in\bbD
\end{alignat}

%%%%%%%%%%%%%%%%%%%%%%%%%%%%%%%%%%%%%%%%%%%%%%%%%%%%%%%
\section{Verblunsky's and Geronimus' Theorems} \lb{s3} 
%%%%%%%%%%%%%%%%%%%%%%%%%%%%%%%%%%%%%%%%%%%%%%%%%%%%%%%

In this section, we will prove Verblunsky's theorem (Theorem~\ref{T2.2}) and 
also a celebrated theorem of Geronimus. Our approach follows Section~3.1 
of \cite{OPUC1} which claims a new proof of Geronimus' theorem assuming 
Verblunsky's theorem. But in preparing this article, we realized the argument 
can be slightly modified to also prove Verblunsky's theorem. 

To state Geronimus' theorem, we need to describe the Schur algorithm 
\cite{Schur}. Given a Schur function, $f$, define  
\begin{equation} \lb{3.1} 
\gamma_0 (f) =f(0) \qquad f(z) = 
\f{\gamma_0 + zf_1(z)}{1+\bar\gamma_0 zf_1 (z)} 
\end{equation} 
If $\gamma_0\in\partial\bbD$ (i.e., $f(z)\equiv \gamma_0)$, we do not define 
$f_1$. Otherwise, $f_1$ defined by \eqref{3.1} is also a Schur function since 
$w\to (\gamma_0 +w)/(1+\bar\gamma_0 w)$ is a biholomorphic bijection of 
$\bbD$ to $\bbD$ if $\abs{\gamma_0}<1$, and $g$ a Schur function with $g(0) =0$ 
implies $g(z)/z$ is a Schur function (the Schwarz lemma). 

\eqref{3.1} is called the Schur algorithm. It can be iterated, that is, 
we define $\gamma_n(f)$, the {\it Schur parameters}, and $f_{n+1}$, the 
{\it Schur iterates}, inductively by 
\begin{equation} \lb{3.2} 
\gamma_n(f)=f_n(0) \qquad f_n(z) = 
\f{\gamma_n + zf_{n+1}(z)}{1+\bar\gamma_n zf_{n+1}(z)} 
\end{equation} 
If, for some $n$, $f_n(z)=e^{i\theta_0}$, we set $\gamma_n =e^{i\theta_0}$ 
and stop. In this way, we map any Schur function, $f$, to a sequence in 
$\bbD^{\infty,c}$. We can now state Geronimus' theorem: 

\begin{theorem}[Geronimus' Theorem] \lb{T3.1} Let $\mu$ be a probability 
measure on $\partial\bbD$, $f$ its Schur function, and $\gamma_n(d\mu)
\equiv\gamma_n(f)$ the Schur parameters of $f$. Then 
\begin{equation} \lb{3.2a} 
\gamma_n (d\mu) =\alpha_n (d\mu) 
\end{equation} 
\end{theorem} 

This gives a continued fraction expansion of $F$ whose coefficients are $\alpha_n$, 
and so is an analog of \eqref{1.6}. This formula explains why we took a minus and 
conjugate in \eqref{2.16}. The procedure of dropping a Verblunsky coefficient 
from the start can be understood by using the recursion relations and the 
relation of $F$ to the OPUC (see Theorem~\ref{T4.5} below). This approach to 
proving Theorem~\ref{T3.1}, due to Peherstorfer \cite{Pe96}, is discussed 
in Section~3.3 of \cite{OPUC1}. 

\eqref{3.1}/\eqref{3.2} can be rewritten and then iterated following 
Schur \cite{Schur}: 
\begin{align} 
f(z) &= \gamma_0 + (1-\bar\gamma_0 f)zf_1 \lb{3.3} \\
&= \gamma_0 + \sum_{j=1}^{n-1} \biggl[\, \prod_{k=0}^{j-1} 
(1-\bar\gamma_k f_k)\biggr] z^j \gamma_j + \prod_{k=0}^{n-1} 
(1-\bar\gamma_k f_k)z^n f_n \lb{3.4}  
\end{align}
which implies that if $f(z) =\sum_{n=0}^\infty a_n (f)z^n$, then 
\begin{equation} \lb{3.5} 
a_n(f) =\gamma_n \prod_{j=0}^{n-1} (1-\abs{\gamma_j}^2) + 
\text{ polynomial in } (\gamma_0, \bar\gamma_0, \dots, 
\gamma_{n-1}, \bar\gamma_{n-1}) 
\end{equation} 
Plugging this into \eqref{2.4} and using \eqref{2.6} implies 
\begin{equation} \lb{3.6} 
c_n(d\mu) =\gamma_{n-1} \prod_{j=0}^{n-2} (1-\abs{\gamma_j}^2) 
+ \text{ polynomial in } (\gamma_0, \bar\gamma_0, \dots, \gamma_{n-2}, 
\bar\gamma_{n-2})
\end{equation}
(the polynomials are different but the leading terms are the same up to 
a shift of index). 

\eqref{3.5} also shows that if $\gamma_j(f)=\gamma_j(g)$ for $j=0, \dots, 
n-1$, then the Schur function $\f12 (f-g)=O(z^n)$ so, by the Schwarz lemma, 
\begin{equation} \lb{3.7} 
\gamma_j (f) =\gamma_j(g), \quad j=0,\dots, n-1 \Rightarrow 
\abs{f(z)-g(z)} \leq 2\abs{z}^n 
\end{equation}

\begin{lemma}\lb{L3.2} The map from Schur functions to $\bbD^{\infty,c}$ is 
one-one and onto. 
\end{lemma} 

\begin{proof} \eqref{3.7} shows that if $\gamma_j(f)=\gamma_j(g)$ for all 
$j$, then $f=g$ on $\bbD$. Given a sequence in $\bbD^\infty$, define 
the {\it Schur approximates}, $f^{[n]}$, by setting $f_{n+1}^{[n]}$ to $0$ 
in \eqref{3.2} and using $\{\gamma_j\}_{j=0}^n$ to define $f_n^{[n]}, 
f_{n-1}^{[n]}, \dots, f^{[n]}$. By construction, 
\begin{equation} \lb{3.8a}
\gamma_j (f^{[n]}) = \begin{cases} \gamma_j & j\leq n \\
0 & j>n \end{cases} 
\end{equation}
Since $\gamma_j (f^{[n]})=\gamma_j (f^{[m]})$ for $j\leq \min (n,m)$, we 
have, by \eqref{3.7}, that $f^{[n]}$ converge uniformly on compacts and 
the limit clearly has the prescribed set of $\gamma$'s. Given a sequence 
in $\bbD^{\infty, c}\backslash\bbD^\infty$, suppose $\gamma_{n+1}= 
e^{i\theta_0}\in\partial\bbD$, set $f_{n+1}\equiv e^{i\theta_0}$ and 
use \eqref{3.2} to define $f$ with the prescribed $\gamma$'s. 
\end{proof} 

\begin{proof}[Proof of Theorems~\ref{T2.2} and \ref{T3.1}] $(z^n -\Phi_n) \perp 
\Phi_n\Rightarrow \|\Phi_n\|^2 = \langle z^n,\Phi_n\rangle$, so applying 
$\,^{*,n}$, 
\begin{equation} \lb{3.8} 
\langle \Phi_n^*,1\rangle = \|\Phi_n\|^2 = \prod_{j=0}^{n-1} (1-\abs{\alpha_j}^2) 
\end{equation} 
Taking the inner product of \eqref{2.16} with the function, $1$, and using 
$\langle \Phi_{n+1},1\rangle =0$, we see
\begin{equation} \lb{3.9} 
\langle z\Phi_n, 1\rangle = \alpha_n \prod_{j=0}^{n-1} (1-\abs{\alpha_j}^2) 
\end{equation} 

By \eqref{2.16} and induction, the coefficients of $\Phi_j$ are polynomials
$\alpha_0, \bar\alpha_0, \dots, \alpha_{j-1}, \bar\alpha_{j-1}$ and so, 
by induction, the moments $c_{j+1}$ are polynomials in the same $\alpha$'s. 
Then \eqref{3.9} becomes (a formula of Verblunsky) 
\begin{equation} \lb{3.10} 
c_{n+1} =\alpha_n \prod_{j=0}^{n-1} (1-\abs{\alpha_j}^2) + 
\text{ polynomial in } (\alpha_0, \bar\alpha_0, \dots, \alpha_{n-1}, 
\bar\alpha_{n-1}) 
\end{equation} 

We will now prove Theorem~\ref{T3.1} by induction and then Theorem~\ref{T2.2} 
follows from Lemma~\ref{L3.2}. For $n=0$, we have, by \eqref{3.10} and \eqref{3.6}, 
that 
\begin{equation} \lb{3.11} 
c_1 =\alpha_0 =\gamma_0 
\end{equation} 

Suppose we know $\alpha_j=\gamma_j$ for $j=0,1,\dots, n-1$. We fix those $n$ values 
in $\bbD$ and ask what values of $c_{n+1}$ can occur. By \eqref{3.6}, it is a solid 
disk in $\bbC$ of radius $\prod_{j=0}^{n-1} (1-\abs{\alpha_j}^2)$ since $\gamma_n$ 
can run through $\ol{\bbD}$. The center of the disk is some fixed point (given 
fixed $\{\gamma_j\}_{j=0}^{n-1}$). 

By \eqref{3.10}, it is also a subset of the disk of radius $\prod_{j=0}^{n-1} 
(1-\abs{\alpha_j}^2)$ with possibly another center. But since the sets are the 
same, the centers must be the same, and all $\alpha_j$ must occur. Once we 
know the centers and radii are the same, the equality of the formulae for 
$c_{n+1}$ implies $\alpha_n =\gamma_n$. 
\end{proof}

%%%%%%%%%%%%%%%%%%%%%%%%%%%%%%%%%%%%%%%%%%%%%%%%%%%%%%%
\section{Zeros, the Bernstein-Szeg\H{o} Approximation,\\
and Boundary Conditions} \lb{s4} 
%%%%%%%%%%%%%%%%%%%%%%%%%%%%%%%%%%%%%%%%%%%%%%%%%%%%%%%

Our first goal in this section is to prove that the zeros of OPUC lie in 
$\bbD$. There are six proofs of this in \cite{OPUC1}. We pick the one that 
is shortest, using the same argument that led to \eqref{2.17}. 

\begin{theorem}\lb{T4.1} $\Phi_n$ has all its zeros in $\bbD$ and $\Phi_n^*$ 
has all its zeros in $\bbC\backslash\ol{\bbD}$. 
\end{theorem} 

\begin{proof} (Landau \cite{Land}) Let $\Phi_n(z_0)=0$ and define $P(z) = 
\Phi_n(z)/(z-z_0)$. Since $\deg P=n-1$, $P\perp \Phi_n$. Thus 
\begin{equation} \lb{4.1} 
\|P\|^2 = \|zP\|^2 = \|z_0 P+\Phi_n\|^2 = \abs{z_0}^2 \|P\|^2 + \|\Phi_n\|^2 
\end{equation} 
so $\|\Phi_n\|^2 = (1-\abs{z_0}^2)\|P\|^2$, implying $\abs{z_0}<1$. Since 
$\Phi_n^*(z_0) = 0 \Leftrightarrow \Phi_n (1/\bar z) =0$, the result for $\Phi_n$ 
implies the result for $\Phi_n^*$.  
\end{proof} 

Next, we will identify measures with $\alpha_j (d\mu)=0$ for $j\geq n_0$. The 
key is a calculation that goes back to Erd\'elyi et al.\ \cite{ENZG91}. 

\begin{proposition}\lb{P4.3} Let $P_n$ be a polynomial of degree $n$ with all 
zeros in $\bbD$. Let 
\begin{equation} \lb{4.2} 
d\mu = \f{c\, d\theta}{2\pi \abs{P_n(e^{i\theta})}^2}  
\end{equation} 
where $c$ is picked to make $d\mu$ a probability measure. Then for all integral 
$j<n$ {\rm{(}}including $j<0${\rm{)}}, 
\begin{equation} \lb{4.3} 
\langle z^j, P\rangle_{L^2 (\partial\bbD, d\mu)} =0 
\end{equation} 
\end{proposition} 

\begin{proof} 
\begin{align*} 
\langle z^j, P\rangle_{L^2 (\partial\bbD, d\mu)} &= \int e^{-ij\theta} P(e^{i\theta}) \, 
\f{d\theta}{\left. 2\pi z^{-n} P^*(z) P(z)\right|_{z=e^{i\theta}}} \\
&= \f{1}{2\pi i} \oint z^{n-j-1}\, \f{dz}{P^*(z)}
\end{align*} 
is zero for $n-j-1 \geq 0$ since $P^*(z)$ is nonvanishing on $\ol{\bbD}$. 
\end{proof} 

\begin{theorem}\lb{T4.4} Let $d\mu$ be a nontrivial probability measure on 
$\partial\bbD$. Let 
\begin{equation} \lb{4.4} 
d\mu_n = \f{d\theta}{2\pi\abs{\varphi_n (e^{i\theta},d\mu)}^2} 
\end{equation} 
Then $d\mu_n$ is a probability measure with 
\begin{equation} \lb{4.5} 
\alpha_j (d\mu_n) = \begin{cases} \alpha_j (d\mu) & j\leq n-1 \\
0 & j\geq n \end{cases} 
\end{equation} 
\end{theorem} 

\begin{proof} Let $d\nu =c\, d\mu_n$ where $c$ is picked so that $\int d\nu =1$ 
(eventually, we will prove $c=1$). By Proposition~\ref{P4.3}, $\langle z^j, 
\Phi_n (\dott;d\mu)\rangle_{L^2 (\partial\bbD, d\nu)}=0$ for $j=0,1,\dots, n-1$, 
so $\Phi_n (z;d\nu) =\Phi_n (z;d\mu)$. It follows from Theorem~\ref{2.3} that 
$\alpha_j (d\nu)=\alpha_j(d\mu)$ for $j=0, \dots, n-1$ and $\varphi_n (z;d\mu) 
=\varphi_n (z;d\nu)$. Therefore, $1=\int\abs{\varphi_n}^2 \, d\nu = c$, so 
$d\nu =d\mu_n$. 

By Proposition~\ref{4.3}, for any $k\geq 0$, 
\begin{equation} \lb{4.6} 
\langle z^j, z^k \Phi_n \rangle_{L^2 (\partial\bbD, d\mu_n)} =0 
\qquad j=0,\dots, n+k-1  
\end{equation} 
It follows that $\Phi_{n+k} (z;d\mu_n) =z^k \Phi_n (z;d\mu_n)$ and thus, 
$\Phi_{n+k}(0)=0$ for $k\geq 1$. Therefore, by \eqref{2.16a}, $\alpha_j 
(d\mu_n) =0$ for $j\geq n$. 
\end{proof}  

Even though Theorem~\ref{T4.4} was proven by Verblunsky \cite{V36} and rediscovered 
by Geronimus \cite{Ger46} (to whom it is often credited), $d\mu_n$ are called 
Bernstein-Szeg\H{o} approximations since Szeg\H{o} \cite{Sz19} first considered 
measures of this form \eqref{3.2} and Bernstein \cite{Bern2} their OPRL analog. 
Since, for each fixed $j$, $\alpha_j(d\mu_n)\to \alpha_j (d\mu)$ (indeed, they 
are equal for $n>j$), $d\mu_n\to d\mu$ weakly since $\calS$ is a homeomorphism. 

Some thought about the form of $d\mu_n$ suggests its Carath\'eodory function 
should be a rational function whose denominator is $\varphi_n^*$. We will 
prove this by identifying the numerator. The {\it second kind polynomials}, 
$\psi_n$, are the OPUC for the measure $d\mu_{-1}$ with $\alpha_j (d\mu_{-1}) 
=-\alpha_j (d\mu)$. Notice that in terms of the matrix $A$ of \eqref{2.20}, 
\begin{equation} \lb{4.7} 
\binom{\psi_{n+1}}{-\psi_{n+1}^*} = A(z, \alpha_n(d\mu))  
\binom{\psi_n}{-\psi_n^*}
\end{equation} 
(note $\alpha_n (d\mu)$, not $\alpha_n (d\mu_{-1})$). Thus 
\begin{equation} \lb{4.8} 
\left(\begin{array}{rr} \psi_n & \varphi_n \\ - \psi_n^* & \varphi_n^*\end{array} 
\right) =A(z,\alpha_{n-1}) \dots A(z,\alpha_0) \left(\begin{array}{rr} 
1 & 1 \\ -1 & 1 \end{array}\right)
\end{equation}  
Taking determinants, using $\det(A)=z$, 
\begin{equation} \lb{4.9} 
\varphi_n^*\psi_n + \varphi_n \psi_n^* = 2z^n 
\end{equation} 

\begin{theorem}[Verblunsky \cite{V36}]\lb{T4.5} Let $d\mu_n$ be given by 
\eqref{4.4}. Then 
\begin{equation} \lb{4.10} 
F(z,d\mu_n) = \f{\psi_n^*(z;d\mu)}{\varphi_n^*(z;d\mu)} 
\end{equation} 
\end{theorem} 

\begin{proof} For $z=e^{i\theta}$, \eqref{4.9} can be rewritten as 
$\Real (\,\ol{\varphi_n (e^{i\theta})}\, \psi_n (e^{i\theta})) =1$. Thus, 
if $G(z)$ is the right side of \eqref{4.10}, 
\begin{equation} \lb{4.11} 
\Real (G(e^{i\theta})) = \f{1}{\abs{\varphi_n (e^{i\theta})}^2}
\end{equation} 
Since $G$ is analytic in a neighborhood of $\ol{\bbD}$, $\Real G >0$ on 
$\bbD$. Since $G(0)=1$, the complex Poisson representation (see Rudin 
\cite[pg.~235]{Rudin}) and \eqref{4.10} imply that $G(z)$ is the Carath\'eodory 
function of $d\mu_n$. 
\end{proof} 

It is useful to think of $d\mu$ and $d\mu_{-1}$ as embedded in a family 
$d\mu_\lambda$ for $\lambda\in\partial\bbD$. The {\it Aleksandrov family} 
associated to $d\mu$ is defined by 
\begin{equation} \lb{4.12} 
\alpha_j (d\mu_\lambda) = \lambda \alpha_j (d\mu) 
\end{equation} 
Given Geronimus' theorem (Theorem~\ref{T3.1}), it is easy to see that 
\begin{equation} \lb{4.13} 
f(z,d\mu_\lambda) =\lambda f(z,d\mu) 
\end{equation} 
(for $\gamma_0 (\lambda f) =\lambda \gamma_0(f)$ and $(\lambda f)_1 
=\lambda(f_1)$). So, by \eqref{2.4} and its inverse, $zf(z)=(F(z)-1)/
(F(z)+1)$, 
\begin{equation} \lb{4.14} 
F(z,d\mu_\lambda) =\f{(1-\lambda) + (1+\lambda) F(z,d\mu)}
{(1+\lambda) + (1-\lambda) F(z,d\mu)}
\end{equation} 
which is the original definition of Aleksandrov \cite{Alex}; it is Golinskii-Nevai 
\cite{GN}  who realized its relevance to OPUC and boundary conditions. 
If $\varphi_n^{(\lambda)} (z)=\varphi_n (z;d\mu_\lambda)$, then 
\begin{equation} \lb{4.15} 
\binom{\varphi_{n+1}^{(\lambda)}}{\bar\lambda \varphi_{n+1}^{(\lambda)}} = 
A(z,\alpha_n) \binom{\varphi_n^{(\lambda)}}{\bar\lambda \varphi_n^{(\lambda)}} 
\end{equation} 
so $\varphi_n$ and $\varphi_n^{(\lambda)}$ obey the same difference equation, 
but the $n=0$ boundary values change from $\binom{1}{1}$ to $\binom{1}{\bar\lambda}$. 
The Aleksandrov family is the analog of variation of boundary conditions in 
second-order ODE's. 

A direct calculation (via contour integrals) shows that if $\Real (a) >0$, then 
\begin{equation} \lb{4.16} 
\int_0^{2\pi} \f{(1-e^{i\theta}) + (1+e^{i\theta})a} 
{(1+e^{i\theta}) + (1-e^{i\theta})a} \, \f{d\theta}{2\pi} =1 
\end{equation} 
Since $1$ is the Carath\'eodory function of $d\theta/2\pi$, \eqref{4.16} and 
\eqref{4.14} imply 

\begin{theorem}[Aleksandrov \cite{Alex}, Golinskii \cite{Golppt}]\lb{T4.6} 
For the Aleksandrov family, we have 
\begin{equation} \lb{4.17} 
\int_\theta [d\mu_{e^{i\theta}}(\varphi)] \, \f{d\theta}{2\pi} = 
\f{d\varphi}{2\pi} 
\end{equation} 
\end{theorem} 

This is the OPUC analog of the Javrjan \cite{Javr}-Wegner \cite{Wegner} averaging for 
Schr\"odinger operators, which is the basis of the localization proof of Simon-Wolff 
\cite{SimWolff}. It can be used \cite{OPUC2} to prove localization for suitable random OPUC.

%%%%%%%%%%%%%%%%%%%%%%%%%%%%%%%%%%%%%%%%%%%%%%%%%%%%%%%
\section{The CMV Matrix} \lb{s5} 
%%%%%%%%%%%%%%%%%%%%%%%%%%%%%%%%%%%%%%%%%%%%%%%%%%%%%%%

Perturbation theory involves looking at similarities of measures when their Verblunsky 
coefficients are close in some suitable sense. In the analogous OPRL situation, 
the Jacobi matrices, \eqref{1.4}, are an invaluable tool. If one defines the essential 
support of a measure to be the support with isolated points removed, and if $\rho$ and 
$\gamma$ are measures on $[c,d]\subset\bbR$ with Jacobi parameters $a_n,b_n$ and 
$\ti a_n, \ti b_n$, then $\rho $ and $\gamma$ have the same essential support if 
$\abs{a_n -\ti a_n} + \abs{b_n -\ti b_n}\to 0$. This can be seen by noting that 
the difference of the Jacobi matrices is compact and then  appealing to Weyl's 
theorem on the invariance of essential spectrum. 

In this section, we discuss a suitable matrix representation for multiplication 
by $z$ in $L^2 (\partial\bbD, d\mu)$. There is an obvious choice, namely, $\langle 
\varphi_n, z\varphi_m\rangle$, but this is not the ``right" one. It has two problems. 
If $\sum\abs{\alpha_j}^2 <\infty$, $\{\varphi_n\}_{n=0}^\infty$ is not a basis, and 
so this matrix is not unitary. Even worse, this matrix has finite columns ($\langle 
\varphi_n, z\varphi_m\rangle =0$ if $n>m+1$) but, in general, it does not have finite 
rows. 

The right basis, as discovered by Cantero, Moral, and Vel\'azquez \cite{CMV} is 
the one, $\chi_0,\chi_1,\chi_2, \dots$, obtained by orthonormalizing $1,z,z^{-1},
z^2, z^{-2}, \dots$. We will also want to consider the basis, $x_0, x_1, x_2, 
\dots$ obtained by orthonormalizing $1, z^{-1}, z, z^{-2}, \dots$. Remarkably, 
the $\chi$'s can be expressed in terms of $\varphi$'s and $\varphi^*$'s, and 
the matrix elements in terms of $\alpha$'s and $\rho$'s. 

\begin{proposition} \lb{P5.1} 
\begin{alignat}{3} 
&\text{\rm{(a)}}\quad \qquad && \chi_{2n}(z) = z^{-n} \varphi_{2n}^*(z) && 
\qquad \chi_{2n-1}(z) = z^{-n+1} \varphi_{2n-1}(z)\lb{5.1} \\
&\text{\rm{(b)}} \quad \qquad && x_{2n}(z) = z^{-n} \varphi_{2n}(z) && 
\qquad x_{2n-1}(z) = z^{-n} \varphi_{2n-1}^*(z) \lb{5.2}
\end{alignat}
\end{proposition}

\begin{proof} In terms of the projections $P_{[k,\ell]}$ of \eqref{2.18a}, we have 
\begin{equation} \lb{5.3x}
\varphi_m = \f{(1-P_{[0,m-1]})z^m}{\| \ldots \|} \qquad 
\varphi_m^* = \f{(1-P_{[1,m]})1}{\| \ldots \|}
\end{equation} 
where $\|\dots\|$ is the norm of the numerator. Since multiplication by $z^\ell$ 
is unitary, 
\[
z^{-n} \varphi_{2n} = \f{(1-P_{[-n,n-1]})z^n}{\|\dots\|} = x_{2n}
\]
proving the first half of \eqref{5.2}. The others are similar. 
\end{proof} 

We define four matrices $(\calC =$ {\it CMV matrix}) by: 
\begin{equation} \lb{5.3}
\calC_{k\ell} = \langle\chi_k, z\chi_\ell\rangle \quad 
\ti\calC_{k\ell} = \langle x_k, zx_\ell\rangle \quad 
\calL_{k\ell} =\langle \chi_k, zx_\ell\rangle \quad 
\calM_{k\ell} =\langle x_k, \chi_\ell\rangle
\end{equation} 
Clearly, 
\begin{equation} \lb{5.4}
\calC =\calL\calM \qquad \ti\calC=\calM\calL \qquad 
\ti\calC = \calC^t
\end{equation} 
where the last comes from the fact that the explicit formulae below show $\calL$ 
and $\calM$ are (complex) symmetric. Define, for $\alpha\in\ol{\bbD}$, the 
$2\times 2$ symmetric matrix: 
\begin{equation} \lb{5.5x}
\Theta(\alpha) = \left(\begin{array}{rr} 
\bar\alpha & \rho \\ \rho & -\alpha \end{array}\right)
\end{equation} 

\begin{theorem}\lb{T5.2} Let $\boldsymbol{1}$ be the $1\times 1$ unit matrix. Then 
\begin{equation} \lb{5.4x}
\calM = \boldsymbol{1} \oplus \Theta (\alpha_1)\oplus\Theta(\alpha_3)\oplus\cdots 
\quad 
\calL = \Theta(\alpha_0)\oplus\Theta (\alpha_2)\oplus\Theta (\alpha_4)\oplus \cdots
\end{equation} 
\end{theorem}

\begin{proof} This is an expression of the Szeg\H{o} recursion formula. For example, 
the $2n$ row (labelling rows $0,1,2,\dots$) of $\calL$ says that $zx_{2n} = 
\bar\alpha_{2n} \chi_{2n} + \rho_{2n} \chi_{2n+1}$ which, by Proposition~\ref{P5.1}, 
is equivalent to $z\varphi_{2n}=\bar\alpha_{2n} \varphi_{2n}^* + \rho_{2n} 
\varphi_{2n+1}$, which is the top row of \eqref{2.19b}. 
\end{proof} 

While $\calL$ and $\calM$ have direct sum structures, in general (i.e., if all 
$\abs{\alpha_j} <1$), $\calC$ does not. Indeed, by \eqref{5.4} and \eqref{5.4x}, 
\begin{equation} \lb{5.5}
\calC = \begin{pmatrix} 
{}& \bar\alpha_0 & \bar\alpha_1 \rho_0 & \rho_1 \rho_0 & 0 & 0 & \dots & {} \\
{}& \rho_0 & -\bar\alpha_1 \alpha_0 & -\rho_1 \alpha_0 & 0 & 0 & \dots & {} \\
{}& 0 & \bar\alpha_2\rho_1 & -\bar\alpha_2 \alpha_1 & \bar\alpha_3 \rho_2 & \rho_3 \rho_2 & \dots & {} \\
{}& 0 & \rho_2 \rho_1 & -\rho_2 \alpha_1 & -\bar\alpha_3 \alpha_2 & -\rho_3 \alpha_2 & \dots & {} \\
{}& 0 & 0 & 0 & \bar\alpha_4 \rho_3 & -\bar\alpha_4 \alpha_3 & \dots & {} \\ 
{}& \dots & \dots & \dots & \dots & \dots & \dots & {} 
\end{pmatrix}
\end{equation}
Thus $\calC$ has a $4\times 2$ block structure and is generally five-diagonal. 
It is the simplest unitary matrix with a cyclic vector; for example \cite{BHJ}, 
any tridiagonal semi-infinite unitary is a direct sum of $1\times 1$ and 
$2\times 2$ matrices. 

\begin{theorem}\lb{T5.3} If $\calC^{(N)}$ is the top left $N\times N$ block of $\calC$, 
then 
\begin{equation} \lb{5.6}
\det(z_0\boldsymbol{1} -\calC^{(N)}) = \Phi_N (z_0)
\end{equation} 
\end{theorem} 

\begin{proof}[Sketch] If $\zeta$ is the operator of multiplication by $z$ in 
$L^2 (\partial\bbD, d\mu)$, then $\calC^{(N)}=P_{[-\ell, N-\ell-1]} \zeta 
P_{[-\ell, N-\ell-1]}$ restricted to $\ran P_{[-\ell, N-\ell-1]}$ where 
$P_{[j,k]}$ is given by \eqref{2.18a} and $\ell$ is either $(N-1)/2$ or 
$N/2$. Since multiplication by $z^\ell$ is unitary, $\calC^{(N)}$ is unitarily 
equivalent to $P_{[0,N-1]}\zeta P_{[0,N-1]}$ on $\ran P_{[0,N-1]}$. 

$z_0$ is an eigenvalue of $P_{[0,N-1]}\zeta P_{[0,N-1]}$ if and only if 
there is $Q$ of degree $N-1$ so $(z-z_0)Q=\Phi_N(z)$, that is, if and only if 
$\Phi_N(z_0)=0$. This proves \eqref{5.6} if $\Phi_n$ has distinct zeros. By a 
limiting argument (see Theorem~1.7.18 of \cite{OPUC1}), \eqref{5.6} holds in 
general. 
\end{proof} 

{\it Remark.} This theorem sheds light on a result of Fej\'er \cite{Fej} that 
for OP's of general measures on $\bbC$, their zeros lie in the convex hull of 
$\supp (d\mu)$. For \eqref{5.6} implies the zeros are in the numerical range of 
$\calC^{(N)}$, so the numerical range of $\calC$, which is the convex hull 
of $\supp (d\mu)$ by the spectral theorem. In particular, Fej\'er's theorem 
implies in the OPUC case that if $\zeta\in\partial\bbD$ with $d=\dist (\zeta, 
\supp(d\mu))>0$ and $\Phi_n (z_0)=0$, then $\abs{z_0-\zeta}\geq \f12 d^2$. 

Notice that if $\abs{\alpha_j}=1$, $\Theta(\alpha_j)= \left( \begin{smallmatrix} 
\bar\alpha_j & 0 \\ 0 & -\alpha_j \end{smallmatrix}\right)$ is a direct sum, and so 
$\calC =\calL\calM$ has a $(j+1)\times (j+1)$ unitary block in the upper corner. 
This implies that if $\beta\in\partial\bbD$, then $\Phi_N^{(\beta)}\equiv z 
\Phi_{N-1}-\beta\Phi_{N-1}^*$ has all its zeros on $\partial\bbD$ since they are 
eigenvalues of a unitary matrix. The $\Phi_N^{(\beta)}$ are called {\it paraorthogonal 
polynomials} and studied in \cite{Ger46,JNT89,Gol02}. 

Dombrowski \cite{Dom78} proved that a Jacobi matrix with $\liminf a_n=0$ has no 
a.c.\ spectrum by picking a subsequence with $\sum_{j=0}^\infty a_{n(j)}<\infty$ 
and trace class perturbing to a decoupled direct sum of finite rank matrices. 
Unaware of this work, Simon-Spencer \cite{S208} proved a similar result if 
$\limsup \abs{b_n} =\infty$. As noted by Golinskii-Simon \cite{GolSim}, the same 
idea and CMV matrices prove the following, originally proven by other means 
\cite{Rakh83}. 

\begin{theorem}[Rakhmanov's Lemma \cite{Rakh83}] \lb{T5.4} If $\mu$ is a 
probability measure on $\partial\bbD$ so $\limsup \abs{\alpha_n(d\mu)}=1$, 
then $\mu$ is singular with respect to $d\theta/2\pi$. 
\end{theorem} 

Golinskii-Simon also use perturbations of CMV matrices to prove: 

\begin{theorem}[\cite{GolSim}] \lb{T5.5} If $\mu,\nu$ are two probability measures on 
$\partial\bbD$ so $\abs{\alpha_n(d\mu) -\alpha_n (d\nu)}\to 0$, then $\esssup 
(d\mu) =\esssup (d\nu)$. If $\sum_n \abs{\alpha_n (d\mu) -\alpha_n (d\nu)}<\infty$, 
then the absolutely continuous parts of $\mu$ and $\nu$ are mutually absolutely 
continuous. 
\end{theorem} 

Aleksandrov families fit into CMV matrices with a twist. $\calC(\{\lambda\alpha_n\})$ 
and $\calC (\{\alpha_n\})$ do not differ by a rank one perturbation --- rather they 
do up to a unitary equivalence. Specifically: 

\begin{theorem} \lb{T5.6} Let $\lambda\in\partial\bbD$ and $\{\alpha_n\}\in\bbD^\infty$. 
Let $D$ be the diagonal matrix with elements $1,\lambda^{-1}, 1, \lambda^{-1}, \dots$. 
Then $D\calC(\{\lambda\alpha_n\})D^{-1} =\calL(\{\alpha_n\})\calM_\lambda (\{\alpha_n\})$ 
where $\calM_\lambda$ differs from $\calM$ by having $\lambda$ in the $(1,1)$ position 
instead of $1$. 
\end{theorem} 

This is a restatement of Theorem~4.2.9 of \cite{OPUC1}. A generalization to rank 
one perturbation in the $n$-th diagonal can be found in Simon \cite{S297}. 

\smallskip
CMV matrices have been generalized in two directions. First, OPUC can be thought 
of as an analog of half-line ODE. The whole-line analog is an extended CMV matrix, 
$\calE$, defined on $\ell_2 (-\infty,\infty)$ by a two-sided sequence 
$\{\alpha_n\}_{n=-\infty}^\infty$ as a product of $\cdots\oplus\Theta (\alpha_{-2}) 
\oplus\Theta (\alpha_0)\oplus\Theta (\alpha_2)\oplus\cdots$ and $\cdots\oplus\Theta 
(\alpha_{-1})\oplus\Theta (\alpha_1) \oplus\cdots$ where $\Theta (\alpha_j)$ acts 
on the span of $\delta_j$ and $\delta_{j+1}$. This is discussed in Sections~4.5, 
10.5, and 10.16 of \cite{OPUC1,OPUC2}. It is useful in the study of ergodic 
(Section~\ref{s6}) and periodic (Section~\ref{s10}) OPUC. Gesztesy-Zinchenko 
\cite{GZ1,GZ2} have further results on $\calE$. 

Second, if $U$ is an $n\times n$ unitary matrix and $\varphi$ is cyclic in that 
$\{U^j\varphi\}_{j=0}^{n-1}$ is a basis, then the spectral measure for $\varphi$ 
has $n$ points, defines polynomials $\Phi_0, \dots, \Phi_n$ and Verblunsky 
coefficients $\alpha_0, \dots, \alpha_{n-2}\in\bbD$ and $\alpha_{n-1}\in\partial\bbD$. 
$U$ is unitarily equivalent to a finite CMV matrix, the upper block of an infinite 
matrix where $\alpha_{n-1}$ is taken in $\partial\bbD$. 

Just as the theory of selfadjoint matrices with cyclic vector is identical to 
the theory of Jacobi matrices, the theory of unitary matrices with cyclic vector 
(i.e., $\{U^j\varphi\}_{j=-\infty}^\infty$ spanning) is identical to the theory 
of CMV matrices. The Verblunsky coefficients are a complete set of unitary 
invariants. 

In this regard, there is a natural question answered by Killip-Nenciu \cite{KilNen}. 
Let $\bbU(n)$ be the group of $n\times n$ unitary matrices and consider Haar measure 
on $\bbU(n)$. For a.e.\ $U$\!, $(1\, 0 \dots 0)^t$ is cyclic, so there is induced a 
measure on Verblunsky coefficients $\alpha_0, \dots, \alpha_{n-2}\in\bbD$ and 
$\alpha_{n-1}\in\partial\bbD$. The measure is the same if $(1\, 0 \dots 0)^t$ is 
replaced by any other vector or by a random choice (say, uniform distribution on 
the unit sphere in $\bbC^n$). 

\begin{theorem}[\cite{KilNen}]\lb{T5.6A} Under the measure induced by Haar measure 
on $\bbU(n)$, the $\alpha_j$ are independent {\rm{(}}i.e., the induced measure 
is a product measure{\rm{)}}, $\alpha_{n-1}$ is uniformly distributed on 
$\partial\bbD$, and $\alpha_j$, $j=1, \dots, n-2$, is distributed via 
\begin{equation} \lb{5.7}
\f{(n-j-1)}{\pi}\, (1-\abs{\alpha}^2)^{(n-j-2)}\, d^2\alpha
\end{equation} 
\end{theorem}

%%%%%%%%%%%%%%%%%%%%%%%%%%%%%%%%%%%%%%%%%%%%%%%%%%%%%%%%%%%%%%%%%%%%%%%%%%%
\section{Transfer Matrices, Weyl Solutions, and Lyapunov Exponents} \lb{s6}
%%%%%%%%%%%%%%%%%%%%%%%%%%%%%%%%%%%%%%%%%%%%%%%%%%%%%%%%%%%%%%%%%%%%%%%%%%% 

In this section, we present a potpourri of results connected with solutions of 
Szeg\H{o} recursion \eqref{2.19b} where the two components are freed of $u_2^*=u_1$  
--- indeed, we look at solutions for a fixed $z$. Thus, solutions have the form 
\begin{equation} \lb{6.1} 
u(z;n)=T_n(z) u(z;0) \qquad T_n(z) =A(z, \alpha _{n-1})\dots A(z,\alpha_0)
\end{equation} 
with $A$ given by \eqref{2.20}. $T_n$ is called the {\it transfer matrix}. 
By \eqref{4.8}, we have 
\begin{equation} \lb{6.2} 
\begin{aligned}
T_n(z) &= \tfrac12 \begin{pmatrix} 
\varphi_n(z) + \psi_n(z) & \varphi_n(z) -\psi_n(z) \\ 
\varphi_n^*(z) - \psi_n^*(z) & \varphi_n^*(z) + \psi_n^*(z) \end{pmatrix} \\
&= \biggl(\,\prod_{j=0}^{n-1} \rho_j^{-1/2}\biggr)
\begin{pmatrix} 
zB_{n-1}^*(z) & -A_{n-1}^*(z) \\
-zA_{n-1}(z) & B_{n-1}(z) \end{pmatrix}
\end{aligned}
\end{equation} 
where $A_{n-1}$ and $B_{n-1}$ are degree $n-1$ polynomials and the $\,^*$ term is 
$\,^{*,n-1}$. The degree count uses $\varphi_n(0)=-\psi_n(0)$, $\varphi_n^*(0) = 
\psi_n^*(z)$. $A_n$ and $B_n$ are the {\it Wall polynomials} which are related to 
the Schur approximants, $f^{[n]}$, of \eqref{3.8a} by $f^{[n]}(z)=A_n(z)/B_n(z)$.

For $z\in\partial\bbD$, $T_n$ lies in the group $\bbU(1,1)$ of matrices obeying 
$M^* \left( \begin{smallmatrix} 1 & 0 \\ 0 & -1\end{smallmatrix}\right)M =
\left( \begin{smallmatrix} 1 & 0 \\ 0 & -1\end{smallmatrix}\right)$. Features of 
this group play a role in advanced aspects of the theory; see \cite{OPUC2}, 
especially Section~10.4. 

The solutions $u_\varphi = (\varphi_n, \varphi_n^*)$ and $u_\psi = (\psi_n, 
-\psi_n^*)$ of \eqref{6.1} can be combined into an $\ell_2$ solution for $\abs{z} 
<1$: 

\begin{theorem}[Geronimo \cite{Ge92}; Golinskii-Nevai \cite{GN}] \lb{T6.1} 
Fix $z\in\bbD$. Then $u_\psi(z;n) + ru_\varphi(z;n)\to 0$ as $n\to\infty$ for 
fixed $r\in\bbC$ if and only if $r=F(z)$. Moreover, $u_\psi + F(z)u_\varphi$ is 
in $\ell^2$. 
\end{theorem} 

{\it Remark.} In analogy to ODE theory, $u_\psi + F(z)u_\varphi$ is called the 
{\it Weyl solution}.

\begin{proof}[Sketch (\cite{GN}) ]Looking at the second component, we see that 
if $u_\psi + ru_\varphi\to 0$, then $-\psi_n^* + r\varphi_n^*\to 0$. By \eqref{2.22}, 
$r-\psi_n^*/\varphi_n^*\to 0$, so by  \eqref{4.10}, $r=F(z)$. That the first 
components go to zero for $r=F(z)$ will follow from the $\ell^2$ proof.  

By using the CD formula \eqref{2.21x} for $\varphi$ and $\psi$ plus a mixed CD 
formula obtained from \eqref{2.30a} by using $\binom{1}{1}$ and $\binom{1}{-1}$, 
one finds that 
\begin{equation} \lb{6.3} 
(1-\abs{z}^2) \sum_{j=0}^{n-1} \, \abs{\psi_j(z) + r\varphi_j(z)}^2 = 4\Real (r) + 
\abs{\psi_n^*(z) - r\varphi_n^*}^2 - \abs{\psi_n(z) + r\varphi_n(z)}^2
\end{equation} 
Taking $r=\psi_n^*(z)/\varphi_n^*(z)$, one finds  
\begin{equation} \lb{6.4} 
k\leq n-1 \Rightarrow \sum_{j=0}^k \, \biggl| \psi_j + \f{\psi_n^*}{\varphi_n^*} 
\varphi_j \biggr|^2 \leq 4 \Real \biggl(\f{\psi_n^*}{\varphi_n^*}\biggr) 
\end{equation} 
Taking $n\to\infty$ and then $k\to\infty$ shows 
\begin{equation} \lb{6.5} 
\sum_{j=0}^\infty\, \abs{\psi_j + F\varphi_j}^2 \leq 4\Real (F) 
\end{equation} 
The inequality in \eqref{6.5} plus equality in \eqref{6.3} imply that 
$\abs{\psi_j^* -F\varphi_j^*} \leq \abs{\psi_j +F\varphi_j}$, so \eqref{6.5} 
implies $u_\psi + Fu_\varphi\in\ell^2$. 
\end{proof} 

Another way of proving the $\ell^2$ result, from \cite{Ge92}, is illuminating. 
It starts from a formula which was Geronimus' original definition of the 
second kind polynomials,  
\begin{equation} \lb{6.6} 
\psi_n(z) = \int \f{e^{i\theta}+z}{e^{i\theta}-z}\, [\varphi_n (e^{i\theta}) - 
\varphi_n(z)]\, d\mu(\theta) 
\end{equation} 
This and its image under the map $\,^*\,$ imply 
\begin{equation}\lb{6.7} 
\begin{aligned} 
F(z) \varphi_n(z) + \psi_n(z) &=\int \f{e^{i\theta}+z}{e^{i\theta}-z}\, 
\varphi_n(e^{i\theta})\, d\mu(\theta) \\
F(z) \varphi_n^*(z) -\psi_n^*(z) &= z^n \int \f{e^{i\theta}+z}{e^{i\theta}-z}\, 
\ol{\varphi_n (e^{i\theta})}\, d\mu(\theta)
\end{aligned} 
\end{equation}
Using $\int\abs{\varphi_n}\, d\mu \leq 1$ and $(e^{i\theta}+z)(e^{i\theta}-z)^{-1} 
=1 + \sum_{n=1}^\infty 2(e^{-i\theta}z)^n$, we see the Taylor coefficients of 
each expression in \eqref{6.7} are bounded by $2$. Since $\int e^{-ik\theta} 
\varphi_n (e^{i\theta})\, d\mu(\theta) =0$ for $k=0, \dots, n-1$, we see 
$\abs{F\varphi_n + \psi_n}\leq 2 \abs{z}^n (1-\abs{z})^{-1}$, while 
$\abs{F\varphi_n^* -\psi_n^*}\leq 2 \abs{z}^{n+1} (1-\abs{z})^{-1}$. This 
proves not only an $\ell^2$ property but exponential decay. 

\smallskip
The next issue we want to discuss is Lyapunov exponents. To understand them, it 
pays to also discuss the {\it density of zeros}, an object of independent interest. 
Given $d\mu$, a nontrivial probability measure in $\partial\bbD$, define the 
measure $d\nu_n$ on $\ol{\bbD}$ to be the point measure which gives weight $k/n$ 
to a zero of $\Phi_n$ of multiplicity $k$. On account of \eqref{5.6} for 
$\ell=0,1,2,\dots$, 
\begin{equation} \lb{6.8} 
\int z^\ell \, d\nu_n(z) = \f{1}{n}\, \tr ([\calC^{(n)}]^\ell) 
\end{equation} 
which can help show that $d\nu_n$ sometimes has a weak limit; if it does, we say 
the limit is the {\it density of zeros}. The limit may not exist; there even 
exist examples (see Example~1.1.17 of \cite{OPUC1}) where the set of limit points 
of $d\nu_n$ is all measures on $\ol{\bbD}$! Here is how \eqref{6.8} can be used: 

\begin{theorem}[Mhaskar-Saff \cite{MhS1}]\lb{T6.2} If 
\begin{equation} \lb{6.9} 
\lim_{n\to\infty} \, \abs{\alpha_n}^{1/n} = r \qquad\text{and}\qquad 
\f{1}{n}\, \sum_{j=0}^{n-1}\, \abs{\alpha_j}\to 0 
\end{equation} 
{\rm{(}}automatic if $r<1${\rm{)}}, then $d\nu_n$ converges to the uniform measure 
on the circle of radius $r$. 
\end{theorem} 

\begin{proof}[Sketch] (See Sections~8.1 and 8.2 of \cite{OPUC1} for details.)  
An argument that exploits Theorem~\ref{T9.1} below and the fact that $\abs{\alpha_{n-1}} 
=\abs{\Phi_n(0)}$ is the product of zeros shows that when the first equation in 
\eqref{6.9} holds, then all limits of $d\nu_n$ are concentrated on the circle of 
radius $r$. The second equation and \eqref{6.8} show that any such limit, $d\nu$, 
has $\int z^\ell \, d\nu =\delta_{\ell 0}$ for $\ell\geq 0$. 
\end{proof} 

The other case where we know $\nu_n$ has a limit is ergodic families of Verblunsky 
coefficients. Let $(\Omega,d\beta)$ be a probability measure space, $T:\Omega\to 
\Omega$, an invertible ergodic transformation, and $V:\Omega\to\bbD$. For each 
$\omega\in\Omega$, define a measure $d\mu_\omega$ by 
\begin{equation} \lb{6.10} 
\alpha_j (d\mu_\omega) = V(T^j\omega) 
\end{equation} 
An argument using the ergodic theorem, \eqref{6.8}, and control of $\lim \abs{\alpha_n 
(d\mu_\omega)}^{1/n}$ show that so long as $\int [-\log V(\omega)]\, d\beta(\omega) 
<\infty$, then $d\mu_\omega$ has for a.e.\ $\omega$ a limit supported on $\partial 
\bbD$ and $\omega$-independent. The most important examples of ergodic families are 
random, periodic, almost periodic, and subshifts (see Chapters~10--12 of \cite{OPUC2}). 

Before leaving the subject of zeros, we note: 

\begin{theorem}[Widom \cite{Wi67}] \lb{T6.3} If $\supp(d\mu)$ is not all of 
$\partial\bbD$, then for any $r<1$, $\sup_n (\#\text{ of zeros of $\Phi_n$ in } 
\abs{z}<r)<\infty$. In particular, any limit of $d\nu_n$ is supported on $\partial\bbD$. 
\end{theorem} 

\begin{theorem}[see Theorem~8.1.11 of \cite{OPUC1}] \lb{T6.4} If $z_0$ in $\partial\bbD$ 
is an isolated point of $\supp(d\mu)$, there is precisely a single zero of $\Phi_n$ near $z_0$ 
for $n$ large and it approaches $z_0$ exponentially fast. 
\end{theorem} 

Finally, we discuss the Lyapunov exponent and Thouless formula: 

\begin{theorem}[see Theorem~10.5.8 of \cite{OPUC2}]\lb{T6.5} If the density of zeros 
measure, $d\nu$, exists and is supported on $\partial\bbD$, and if 
\begin{equation} \lb{6.11} 
\rho_\infty = \lim_{n\to\infty}\, \biggl(\, \prod_{j=0}^{n-1} \rho_j\biggr)^{1/n}  
\end{equation} 
exists, then for $z\notin\partial\bbD$, the following limit exists and is given by 
\begin{equation} \lb{6.12} 
\gamma(z) \equiv\lim_{n\to\infty}\, \f{1}{n}\, \log \|T_n(z)\| =-\log \rho_\infty 
-\int \log \abs{e^{i\theta}-z}^{-1}\, d\nu(\theta)
\end{equation} 
\end{theorem} 

$\gamma$ is called the {\it Lyapunov exponent}. \eqref{6.12} is called the 
{\it Thouless formula}. For $\abs{z}>1$, $\abs{\varphi_n} >\abs{\varphi_n^*}$ and 
$\abs{\psi_n} > \abs{\psi_n^*}$, we need only control the growth of $\abs{\varphi_n}$ 
and $\abs{\psi_n}$. By \eqref{6.8}, $\varphi_n$ and $\psi_n$ have the same density 
of zeros. Writing $\varphi_n =\prod_{j=0}^{n-1} \rho_j^{-1} \prod_{\text{zeros}} 
(z-z_\ell)$ easily yields \eqref{6.12}. 

See \cite{OPUC2} for discussion of when \eqref{6.12} holds on $\partial\bbD$ and for 
further study of ergodic OPUC.

%%%%%%%%%%%%%%%%%%%%%%%%%%%%%%%%%%%%%%%%%%%%%%%%%%%%%%%%%%%%%%%%%%%%%%%%%%%%%%%%
\section{Khrushchev's Formula, CMV Resolvents,\\ and Rakhmanov's Theorem} \lb{s7}
%%%%%%%%%%%%%%%%%%%%%%%%%%%%%%%%%%%%%%%%%%%%%%%%%%%%%%%%%%%%%%%%%%%%%%%%%%%%%%%%

In two remarkable papers \cite{Kh2000,Khr}, Khrushchev found deep 
connections between Schur iterates and the structure of OPUC. A key 
input for the theory is: 

\begin{theorem}[Khrushchev's Formula]\lb{T7.1} The Schur function for 
the measure $\abs{\varphi_n (e^{i\theta},d\mu)}^2\, d\mu(\theta)$ is 
given by $b_n(z) f_n(z)$, where $f_n$ is the $n$-th Schur iterate 
{\rm{(}}by Geronimus' theorem, this is the Schur function of the 
measure with Verblunsky coefficients $\{\alpha_{n+j}\}_{j=0}^\infty${\rm{)}} 
and $b_n$ is the Blaschke product, 
\begin{equation} \lb{7.1} 
b_n (z;d\mu) = \f{\varphi_n (z;d\mu)}{\varphi_n^*(z;d\mu)} 
\end{equation} 
\end{theorem}

{\it Remark.} Khrushchev's formula illuminates \eqref{2.8}. In this trivial 
measure case, $\{z_j\}_{j=1}^{n-1}$ are the zeros of $\Phi_{n-1}$ and 
$e^{i\theta_0}$ is the Schur parameter, $\gamma_{n-1}$. 

\smallskip
In terms of the CMV matrix, this gives a formula for $\langle\delta_n, 
(\calC+z)(\calC-z)^{-1} \delta_n\rangle$, and so when $n=m$ for
\begin{equation} \lb{7.2} 
G_{nm}(z) =\langle \delta_n, (\calC-z)^{-1}\delta_m\rangle
\end{equation} 
the analog of the Green's function in ODE's. Half-line Green's functions for 
ODE's have the form $f_- (\min(x,y)) f_+ (\max(x,y))$ where $f_-$ (resp.\ $f_+$) 
obeys boundary conditions at $x=0$ (resp.\ $x=\infty$). There is an analogous 
formula, due to Simon (even if $n\neq m)$, for $G_{nm}$ in terms of the 
OPUC and Weyl solutions. It can be found in Section~4.4 of \cite{OPUC1} 
and generalizes Theorem~\ref{T7.1}. Other proofs of Theorem~\ref{T7.1} appear 
in Theorem~4.5.10 of \cite{OPUC1} and Theorem~9.2.4 of \cite{OPUC2}. 
The most important consequence of Khrushchev's formula is: 

\begin{theorem}[Khrushchev \cite{Kh2000}]\lb{T7.2} The essential support of 
the a.c.\ part of $d\mu$ is all of $\partial\bbD$ if and only if 
\begin{equation} \lb{7.3} 
\lim_{n\to\infty} \, \int_0^{2\pi} \abs{f_n (e^{i\theta},d\mu)}^2\, 
\f{d\theta}{2\pi} =0
\end{equation} 
\end{theorem} 

Since 
\begin{equation} \lb{7.4} 
\int_0^{2\pi} f_n (e^{i\theta}, d\mu) \, \f{d\theta}{2\pi} = f_n(0) =\alpha_n 
\end{equation} 
an immediate corollary is 

\begin{theorem}[Rakhmanov's Theorem]\lb{T7.3} If the essential support of the 
a.c.\ part of $d\mu$ is all of $\partial\bbD$, then 
\begin{equation} \lb{7.5} 
\lim_{n\to\infty}\, \alpha_n =0 
\end{equation} 
\end{theorem} 

This result is originally due to Rakhmanov \cite{Rakh77,Rakh83,Rakh87} with 
important further developments by M\'at\'e-Nevai-Totik \cite{MNT85a,MNT88,Nev89,Nev91}. 
Bello-L\'opez \cite{BHLL} extended this result to arcs, and Denisov \cite{DenPAMS} to 
OPRL. Here are some other important results of Khrushchev's theory: 

\begin{theorem}\lb{T7.3A} 
\[
\wlim \abs{\varphi_n (e^{i\theta})}^2\, d\mu = \f{d\theta}{2\pi} \Leftrightarrow 
(\forall j) \lim_{n\to\infty}\, \alpha_{n+j} \alpha_n =0 
\] 
\end{theorem} 

{\it Remark.} \eqref{6.8} can be reinterpreted as saying weak Ces\`aro limits 
of $\abs{\varphi_n}^2\, d\mu$ are the density of zeros when the latter is supported 
on $\partial\bbD$; see Section~8.2 of \cite{OPUC1}. 

\begin{theorem}\lb{T7.4} Let $f^{[n]}$ be the Schur approximates {\rm{(}}given 
by \eqref{3.8a}{\rm{)}}. Then 
\begin{equation} \lb{7.6} 
\int \abs{f^{[n]} (e^{i\theta}) - f(e^{i\theta})}^2 \, \f{d\theta}{2\pi} \to 0 
\end{equation}
if and only if either  
\begin{SL} 
\item[{\rm{(i)}}] $d\mu_\ac=0$, that is, $\mu$ is purely singular, or 
\item[{\rm{(ii)}}] $\alpha_n (d\mu)\to 0$. 
\end{SL} 
Moreover, if $\wlim \abs{\varphi_n (e^{i\theta})}^2\, d\mu  =d\theta/2\pi$, then 
\eqref{7.6} holds. 
\end{theorem} 

As a consequence of these theorems, we get a result for sparse $\alpha$'s: 

\begin{corollary}\lb{C7.5} If $\lim_{n\to\infty} \alpha_{n+j} \alpha_n=0$ for 
all $j$, but $\limsup_n \abs{\alpha_n} \neq 0$, then $\mu$ is purely singular 
continuous. 
\end{corollary} 

\begin{theorem}\lb{T7.6} Suppose that uniformly on compacts of $\partial\bbD$, 
\begin{equation} \lb{7.7} 
\lim_{n\to\infty}\, \f{\Phi_{n+1}^* (z)}{\Phi_n^*(z)} =G(z)
\end{equation} 
then either $G(z)\equiv 1$ or else for some $a\in (0,1]$ and $\lambda\in 
\partial\bbD$, 
\begin{equation} \lb{7.8} 
G(z)=\tfrac12\, \bigl[(1+\lambda z) + \sqrt{(1-\lambda z)^2 + 4a^2\lambda z}\,\bigr]
\end{equation}
\end{theorem} 

Note we have that $G\equiv 1$ if and only if $\lim_{n\to\infty} \alpha_{n+j}\alpha_n 
=0$ for all $j$ and that Barrios-L\'opez have proven that \eqref{7.8} holds if and  
only if $\lim_{n\to\infty}\abs{\alpha_n}=a$ and $\lim_{n\to\infty} \alpha_{n+1} 
\alpha_n^{-1}=\lambda$. 

Khrushchev has also described all possible $d\nu$'s that can occur as $\wlim 
\abs{\varphi_n}^2\, d\mu$ (i.e., for which the limit exists) and when they can 
occur (essentially, asymptotically period $1$ or $2$). The analogs of these 
w-limit and ratio asymptotic results for OPRL were found by Simon \cite{Sim2003}.

%%%%%%%%%%%%%%%%%%%%%%%%%%%%%%%%%%%%%%%%%%%%%%%%%%%%
\section{Szeg\H{o}'s and Baxter's Theorems} \lb{s8}
%%%%%%%%%%%%%%%%%%%%%%%%%%%%%%%%%%%%%%%%%%%%%%%%%%%%

Szeg\H{o}'s theorems may well be the most celebrated in OPUC.  While they 
have expressions purely in terms of OPUC objects, for historical reasons, 
one should state them in terms of Toeplitz determinants, $D_n(d\mu)$. 
This is defined as the determinant of the $(n+1)\times (n+1)$ matrix 
$\{c_{k-\ell}\}_{0\leq k,\ell\leq n}$ with $c$ given by \eqref{2.7}. 
$D_n$ is the Gram determinant of $\{z^k\}_{k=0}^n$ since 
$\langle z^k, z^\ell\rangle_{L^2 (d\mu)} =c_{k-\ell}$. The invariance of 
such determinants under triangular change of basis implies (using also  
\eqref{2.18}) 
\begin{equation} \lb{8.1} 
D_n (d\mu) = \prod_{j=0}^n \|\Phi_j\|^2 =\prod_{j=0}^{n-1} 
(1-\abs{\alpha_j}^2)^{n-j} 
\end{equation} 
which immediately implies 
\begin{align} 
F(d\mu) & \equiv \lim_{n\to\infty} \, D_n(d\mu)^{1/n} = 
\prod_{j=0}^\infty (1-\abs{\alpha_j}^2) = \lim_{n\to\infty} \, 
\|\Phi_n\|^2 \lb{8.2} \\
G(d\mu) &\equiv \lim_{n\to\infty}\, \f{D_n(d\mu)}{F(d\mu)^{n+1}} = 
\prod_{j=0}^\infty (1-\abs{\alpha_j}^2)^{-j-1} \lb{8.3} 
\end{align} 
$F$ is always defined, although it may be $0$. $G$ is defined so long as 
$F>0$, that is, so long as $\sum_{j=0}^\infty \abs{\alpha_j}^2 <\infty$. 
$G$ may be infinite and is finite if and only if $\sum_{j=0}^\infty 
j\abs{\alpha_j}^2 <\infty$. 

Szeg\H{o}'s theorems express $F$ and $G$ in terms of the a.c.\ weight, 
$w$, of $d\mu$:  
\begin{equation} \lb{8.4} 
d\mu = w(\theta)\, \f{d\theta}{2\pi} + d\mu_\s 
\end{equation} 
where $w\in L^1 (\partial\bbD, \f{d\theta}{2\pi})$ and $d\mu_\s$ is singular 
with respect to $d\theta/2\pi$. 

\begin{theorem}[Szeg\H{o}'s Theorem]\lb{T8.1} 
\begin{equation} \lb{8.5} 
F(d\mu) =\prod_{j=0}^\infty (1-\abs{\alpha_j}^2) = \exp \biggl( \int \log 
(w(\theta)) \, \f{d\theta}{2\pi}\biggr) 
\end{equation} 
\end{theorem} 

{\it Remark.} Szeg\H{o} proved this when $d\mu_\s =0$ in 1915; the proof below 
is basically his proof in \cite{Sz20-21}. The result does not depend on $d\mu_\s$ 
--- this was shown first by Verblunsky \cite{V36}. \cite{OPUC1,OPUC2} have  
five proofs of Theorem~\ref{T8.1}

\begin{proof}[Sketch when $d\mu_\s =0$] Since $\Phi_n^*$ is nonvanishing on 
$\bbD$ and $\Phi_n^*(0)=1$, $\int \log \abs{\Phi_n^* (e^{i\theta})} 
\f{d\theta}{2\pi} =1$. Thus, by Jensen's inequality, 
\[
\|\Phi_n\|^2 \equiv \int \abs{\Phi_n^*(e^{i\theta})}^2 w(\theta)\, 
\f{d\theta}{2\pi} \geq \exp \biggl( \int \log (w(\theta))\, \f{d\theta}{2\pi}\biggr)
\]
so $F(d\mu) \geq$ RHS of \eqref{8.5}. On the other hand, since $\Phi_n^*$ is the 
projection of $1$ to the complement of $[z,\dots, z^n]$ in $[1,\dots,z^n]$, 
we have 
\begin{equation} \lb{8.6} 
\|\Phi_n^*\|^2 = \min \{\|P\|_{L^2 (d\mu)}^2 \mid \deg P\leq n,\, P(0)=1\}
\end{equation} 
Using \eqref{8.2} and a limit argument, 
\begin{equation} \lb{8.7} 
F(d\mu) =\min \{\|f\|_{L^2(d\mu)}^2 \mid f\in H^\infty, \, f(0)=1\} 
\end{equation} 
Pick the trial functions $f_\veps(z)=g_\veps(z)/g_\veps (0)$ where 
\[
g_\veps (z) =\exp \biggl( -\int \f{e^{i\theta}+z}{e^{i\theta}-z}\, 
\log ( w(\theta) + \veps)\f{d\theta}{4\pi}\biggr)
\]
and take $\veps\downarrow 0$ to get $F(d\mu) \leq$ RHS of \eqref{8.5}. 
\end{proof} 

Because their singular continuous part is arbitrary, once an $\ell^2$ condition 
is dropped, $d\mu$ can be arbitrarily ``bad": 

\begin{theorem} \lb{T8.2A} Let $d\rho$ be a measure on $\partial\bbD$ with 
support all of $\partial\bbD$. Then there exist $d\mu$, a probability measure 
on $\partial\bbD$ mutually equivalent to $d\rho$, so that for all $p>2$, 
\begin{equation} \lb{8.8a} 
\sum_{n=0}^\infty\, \abs{\alpha_n (d\mu)}^p < \infty
\end{equation}
\end{theorem} 

This is Theorem~2.10.1 of \cite{OPUC1}, proven using ideas of Totik \cite{Totik} 
and the bounds in \eqref{8.5}. 

By \eqref{8.5}, we get one of the gems of spectral theory, equivalences between 
some recursion coefficient property and some spectral measure property: 

\begin{corollary}\lb{C8.2} 
\begin{equation} \lb{8.8} 
\sum_{j=0}^\infty \, \abs{\alpha_j}^2 <\infty \Leftrightarrow \int 
\log(w(\theta))\, \f{d\theta}{2\pi} >-\infty
\end{equation} 
\end{corollary} 

The equivalent conditions \eqref{8.8} are called the {\it Szeg\H{o} condition}. 
When they hold, Szeg\H{o} defined the {\it Szeg\H{o} function} by 
\begin{equation} \lb{8.9} 
D(z) =\exp\biggl( \int \f{e^{i\theta}+z}{e^{i\theta}-z}\, \log (w(\theta))  
\f{d\theta}{4\pi}\biggr) 
\end{equation} 
Standard boundary value theory for the Poisson kernel implies $D(e^{i\theta}) = 
\lim_{r\uparrow 1} D(re^{i\theta})$ exists for $d\theta/2\pi$-a.e.\ $\theta\in 
[0,2\pi)$ and 
\begin{equation} \lb{8.10} 
\abs{D(e^{i\theta})}^2 = w(\theta) 
\end{equation} 

\begin{theorem}[Szeg\H{o} \cite{Sz20-21}]\lb{T8.3} Suppose \eqref{8.8} holds. 
Let $D_\ac (e^{i\theta}) =D(e^{i\theta})$ for a.e.\ $\theta$ and $=0$ on a 
supporting set for $d\mu_\s$. Then  
\begin{alignat}{2}
&\text{\rm{(i)}} \qquad 
&&\int \abs{\varphi_n^*(e^{i\theta}) - D_\ac (e^{i\theta})^{-1}}^2 \, d\mu \to 0  \lb{8.11} \\
&{\rm{(ii)}} \qquad  
&&\int \abs{\varphi_n (e^{i\theta})}^2 \, d\mu_\s \to 0 \lb{8.12} \\
&{\rm{(iii)}} \qquad
&&\varphi_n^*(z) \to D(z)^{-1} \quad \text{uniform on compacts in $\bbD$} \lb{8.13x} 
\end{alignat} 
\end{theorem} 

\begin{proof}[Sketch] A short preliminary argument proves that $D\in\bbH^2 
(\bbD)$. Thus the Cauchy formula holds for $\varphi_n^*D$, so 
\begin{equation} \lb{8.13} 
\int (\varphi_n^* D) (e^{i\theta})\, \f{d\theta}{2\pi} = 
\varphi_n^*(0) D(0)\to 1 
\end{equation} 
since, by \eqref{8.5}, $\varphi_n^*(0) D(0) = \prod_{j=0}^{n-1} (1-\abs{\alpha_j}^2)^{-1/2} 
\prod_{j=0}^\infty (1-\abs{\alpha_j}^2)^{1/2}$. By \eqref{8.10}, 
\[ 
\begin{split}
\int \abs{\varphi_n^* - D(e^{i\theta})^{-1}}^2  & w(\theta)\, \f{d\theta}{2\pi} + 
\int \abs{\varphi_n^* (e^{i\theta})}^2 \, d\mu_\s \\
&\quad = \|\varphi_n^*\|_{L^2(d\mu)}^2 
+ 1-2\Real \text{(LHS of \eqref{8.13})}\to 0 
\end{split}
\]
by \eqref{8.13}. This implies (i) and (ii). This then implies that $D\varphi_n^* \to 
1$ in $L^2 (\partial\bbD, \f{d\theta}{2\pi})$, so by $\bbH^2$ theory, (iii) holds. 
\end{proof} 

\eqref{8.13x} and the related
\begin{equation} \lb{8.14} 
z^{-n} \varphi_n(z) \to \ol{D(1/\bar z)} \text{ on }\bbC\backslash\ol{\bbD}  
\end{equation} 
are called {\it Szeg\H{o} asymptotics}. 

\begin{theorem}[Sharp Form of the Szeg\H{o} Strong Theorem \cite{Sz52,Ib,GoIb}]\lb{T8.4} 
If $d\mu_\s =0$, the Szeg\H{o} condition holds, and $\hat L_k$ are the Fourier coefficients 
of $\log w$, then 
\begin{equation} \lb{8.15} 
G(d\mu) =\prod_{j=0}^\infty (1-\abs{\alpha_j}^2)^{-j-1} =\exp\biggl( \, \sum_{n=0}^\infty 
n\abs{\hat L_n}^2\biggr) 
\end{equation} 
\end{theorem} 

{\it Remark.} Szeg\H{o} \cite{Sz52} proved this when $d\mu$ has certain regularity 
properties. The general result is due to Ibragimov \cite{Ib}; see also \cite{GoIb}. 

Seeing when $G(d\mu) <\infty$ leads to a second gem: 

\begin{corollary}\lb{C8.5} 
\begin{equation} \lb{8.16} 
\sum_{j=0}^\infty j\abs{\alpha_j}^2 <\infty \Leftrightarrow d\mu_\s = 
0 \quad \text{and}\quad \sum_{n=0}^\infty n\abs{\hat L_n}^2 <\infty 
\end{equation} 
\end{corollary} 

This corollary relies also on a theorem of Golinskii-Ibragimov \cite{GoIb} that the 
LHS of \eqref{8.16} $\Rightarrow d\mu_\s =0$. This result plus five distinct proofs 
of Theorem~\ref{T8.4} are found in Chapter~6 of \cite{OPUC1}. A sixth proof is in 
Section~9.10 of \cite{OPUC2}. 

A final gem we want to mention is: 

\begin{theorem}[Baxter's Theorem \cite{Bax,Bax2}]\lb{T8.6} Fix $\ell\geq 0$. 
The following are equivalent: 
\begin{alignat*}{2} 
&\text{\rm{(a)}} \qquad \quad &&\sum_{n=0}^\infty n^\ell \abs{\alpha_n} <\infty \\
&\text{\rm{(b)}} \qquad \quad && d\mu_\s =0, \quad \inf_\theta w(\theta) >0, \quad
\text{and}\quad \sum_{n=0}^\infty n^\ell \abs{c_n} <\infty 
\end{alignat*}
In particular if $d\mu_\s =0$ and $\inf_\theta w(\theta) >0$, then $w$ is $C^\infty$ 
if and only if $\sup_n n^\ell \abs{\alpha_n}<\infty$ for all $\ell\geq 0$. 
\end{theorem} 

This is proven in Chapter~5 of \cite{OPUC1}. 

%%%%%%%%%%%%%%%%%%%%%%%%%%%%%%%%%%%%
\section{Exponential Decay} \lb{s9}
%%%%%%%%%%%%%%%%%%%%%%%%%%%%%%%%%%%%

Suppose for some $R>1$, we have 
\begin{equation} \lb{9.1} 
\abs{\alpha_n} \leq CR^{-n} 
\end{equation} 
By \eqref{2.16} and induction, 
\begin{equation} \lb{9.2} 
\sup_{n,\, \abs{z}=1}\, \abs{\Phi_n(z)} \leq \prod_{j=0}^\infty (1+\abs{\alpha_j})<\infty  
\end{equation} 
so, by the maximum principle and \eqref{2.14}, 
\[
\sup_{\abs{z}\geq 1}\, \abs{z}^{-n} \abs{\Phi_n(z)} = \sup_{\abs{z}\leq 1}\, 
\abs{\Phi_n^*(z)} <\infty 
\]
Thus, by \eqref{2.19a}, if $\abs{z}<R$, 
\begin{equation} \lb{9.3} 
\sum_{n=0}^\infty \, \abs{\Phi_{n+1}^* (z) - \Phi_n^*(z)} \leq \sum_{n=0}^\infty \, 
\abs{\alpha_n} \,  \abs{z}^{n+1} < \infty
\end{equation} 
It follows that $\Phi_n^*(z)$ and so $\varphi_n^*(z)$ converges uniformly on compacts 
of $\{z\mid \abs{z}<R\}$ and so, by \eqref{8.13x}, $D(z)^{-1}$ has a continuation 
to this disk. We thus have one-half of: 

\begin{theorem}[Nevai-Totik \cite{NT89}] \lb{T9.1} Fix $R>1$. The following are 
equivalent:
\begin{SL} 
\item[{\rm{(a)}}] $d\mu_\s =0$, the Szeg\H{o} condition holds, and $D(z)^{-1}$ has  
an analytic continuation to $\{z\mid\abs{z}<R\}$. 
\item[{\rm{(b)}}] $\limsup_{n\to\infty}\, \abs{\alpha_n}^{1/n} = R^{-1}$ 
\end{SL} 
\end{theorem} 

The other direction uses the useful formula, 
\begin{equation} \lb{9.4} 
d\mu_\s =0 \Rightarrow \alpha_n = -D(0)^{-1} \int \, \ol{\Phi_{n+1} (e^{i\theta})}\, 
D(e^{i\theta})^{-1}\, d\mu(\theta) 
\end{equation}

Section~7.1 of \cite{OPUC1} has a complete proof. One can say more (\cite[Section~7.2]{OPUC1}) 
when this holds: $\alpha_n +$ Taylor coefficients of $D(z)^{-1} \, \ol{D(1/\bar z)}$ decays  
as $0(R^{-2n+\veps})$. There is also a lot known about asymptotics of the zeros when there is 
exponential decay (see \cite[Sections~8.1 and 8.2]{OPUC1}, \cite{Saff2,MFMS} and references 
therein).

%%%%%%%%%%%%%%%%%%%%%%%%%%%%%%%%%%%%
\section{Periodic OPUC} \lb{s10}
%%%%%%%%%%%%%%%%%%%%%%%%%%%%%%%%%%%% 

The theory of one-dimensional periodic Schr\"odinger operators (a.k.a.\ Hill's equation) 
and of periodic Jacobi matrices has been extensively developed 
\cite{DubMatNov,FlMcL,Krich1,Krich2,McvM1,vMoer}. In the 1940's, Geronimus \cite{Ger44} 
found the earliest results on OPUC with periodic Verblunsky coefficients, that is, 
for some $p\geq 1$ and $j=0,1,2,\dots$, 
\begin{equation} \lb{10.1} 
\alpha_{j+p} =\alpha_j 
\end{equation}
In particular, the case $\alpha_j\equiv a\in\bbD\backslash\{0\}$ yields OPUC called 
{\it Geronimus polynomials} (see Example~1.6.12 of \cite{OPUC1}). Many of the general  
features for OPUC obeying \eqref{10.1} were found by Peherstorfer and collaborators 
\cite{Pe93,Pe96,Pe01,Pe03,PS6,PS1,PS2,PS3,PS4,PS2000,PS00}. Geronimo-Johnson 
\cite{GJo2,GJo1} have studied almost periodic Verblunsky coefficients. A reworking 
with some new results is Chapter~11 of \cite{OPUC2}, which uses methods mimicking 
the periodic Hill-Jacobi theory. 

We suppose henceforth that $p$ is even. A basic object is the {\it discriminant}, 
\begin{equation} \lb{10.2} 
\Delta(z) =\tr(z^{-p/2} T_p(z)) 
\end{equation} 
where $T_p(z)$ is the transfer matrix given by \eqref{6.2}. The $z^{-p/2}$ is included  
since, by $\det(A)=z$, $\det (z^{-p/2}T_p(z))=1$, and so $z^{-p/2}T_p(z)$ has eigenvalues 
$\f{\Delta}{2}\pm i \sqrt{1-(\f{\Delta}{2})^2}$. In particular, these eigenvalues 
have magnitude $1$, that is, $\sup_m \|T_{mp}(z)\|<\infty$ exactly when $\Delta(z)\in 
[-2,2]$. This is part of: 

\begin{theorem}\lb{T10.1} There exists $\{x_j\}_{j=1}^{2p}$, $\{y_j\}_{j=1}^{2p}$ with 
\[
x_1 < y_1 \leq x_2 < y_2 \leq \cdots \leq x_p < y_p \leq x_1 + 2\pi \equiv x_{p+1}
\]
so that the solutions of $\Delta(z)=2$ {\rm{(}}resp.\ $-2${\rm{)}} are exactly 
$e^{ix_1}, e^{iy_2}, e^{ix_3}, \dots, e^{ix_p}$ {\rm{(}}resp.\ $e^{iy_1}, 
e^{ix_2}, e^{iy_3}, \dots, e^{iy_p}${\rm{)}} and $\Delta(z)\in [-2,2]$ exactly on 
\begin{equation} \lb{10.3} 
B=\bigcup_{j=1}^p \{e^{i\theta} \mid x_j \leq \theta_j \leq y_j\}
\end{equation} 
the {\em bands}. $B$ is the essential support of $d\mu_\ac$ and the only possible 
singular spectrum are mass points which can occur in open gaps {\rm{(}}i.e., 
nonempty sets of the form $\{e^{i\theta}\mid y_j <\theta<x_{j+1}\}${\rm{)}} 
with one {\rm{(}}or zero{\rm{)}} mass point in each gap. 
\end{theorem} 

\begin{theorem}\lb{T10.2} Let $d\rho$ be the equilibrium measure for $B$ {\rm{(}}i.e., 
the minimizer for $\calE(\rho)=\int \log \abs{z-w}^{-1}\, d\rho(z)\, d\rho (w)$ 
with $\supp (d\rho)\subset B$ and $\rho(B)=1${\rm{)}}. Let $C_B$ be the 
logarithmic capacity of $B$ {\rm{(}}i.e., $\exp(-\text{ minimizing value of }
\calE(p))${\rm{)}} and $Q(z)$ the logarithmic potential for $B$ {\rm{(}}i.e., 
$Q(z)=\int \log \abs{z-\omega}^{-1}\, d\rho (\omega)${\rm{)}}. Then 
\begin{SL} 
\item[{\rm{(i)}}] $d\rho$ is the density of zeros for $d\mu$, $-[Q(z) +\log C_B]$ is 
the Lyapunov exponent, and $C_B=\prod_{j=0}^{p-1} (1-\abs{\alpha_j}^2)^{1/2p}$. 
\item[{\rm{(ii)}}] $d\rho$ can be written in terms of the discriminant as 
\begin{equation} \lb{10.4} 
d\rho(\theta) \equiv \f{1}{p} \, \f{\abs{\partial\Delta (e^{i\theta})/\partial\theta}}
{(4-\Delta^2 (e^{i\theta}))^{1/2}}\, \f{d\theta}{2\pi}
\end{equation} 
\item[{\rm{(iii)}}] For each $j=1,2,\dots, p$, 
\begin{equation} \lb{10.5} 
\rho (\{e^{i\theta}\mid x_j\leq \theta\leq y_j\}) =\f{1}{p} 
\end{equation} 
\end{SL}
\end{theorem} 

The proof of this result (see Section~11.1 of \cite{OPUC2}) depends on noting that, 
by the Thouless formula, $\gamma(z)$ is harmonic on $\bbC\backslash B$ and $\gamma(z) 
=0$ on $B$. (iii) is related to half of the following result of Peherstorfer motivated 
by an OPRL result of Aptekarev \cite{Apt}: 

\begin{theorem}[Peherstorfer \cite{Pe93}] \lb{T10.3} Let $B$ be a union of $\ell$ disjoint, 
closed intervals, $B_1, \dots, B_\ell$, in $\partial\bbD$. Then $B$ is the set of bands 
of a period $p$ set of $\alpha$'s if and only if 
\begin{SL} 
\item[{\rm{(1)}}] If $d\rho$ is the equilibrium measure of $B$, then $p\rho(B_j)\in\bbZ$ 
for $j=1, \dots, \ell$. 
\item[{\rm{(2)}}] Let $z_1, z_2, \dots, z_{2p}$ be defined clockwise around the circle 
so that $z_1$ is the lower edge of $B_1$ and the $2p$ points are the $2\ell$ band edges 
and those interior points in a band with $p\rho (B_j)\geq 2$, that divide $B_j$ into 
$p\rho (B_j)$ sets with $\rho$-measure $1/p$, each counted twice. Then $z_1 z_4 z_5  
z_8 z_9 \dots z_{2p} =1$.  
\end{SL} 
\end{theorem} 

If some $\rho (B_j)$ is irrational, then there is no periodic family of $\alpha$'s with 
those bands, but there is an almost periodic set, as proven by Geronimo-Johnson \cite{GJo1} 
(see Section~11.8 of \cite{OPUC2}). 

Given a measure $d\mu$ on $\partial\bbD$ so that \eqref{10.1} holds, the {\it Dirichlet data} 
is defined partly as the $p$ points where $\binom{1}{1}$ is an eigenvector for $T_p(z)$, 
that is, zeros of $\varphi_p^*(z) -\varphi_p (z)$. There is one such point in each gap, 
including closed gaps (i.e., $e^{iy_j}$ when $y_j =x_{j+1}$). If the value is at a gap  
edge, the eigenvalue, $\lambda$, of $z^{-p/2} T_p(z)$ for $\binom{1}{1}$ is $\pm 1$. 
Otherwise, it is in $\bbR\backslash\{0,-1,1\}$. In that case, we add $\sigma_j=\pm 1$ 
to the $j$-th Dirichlet point with $\sigma_j =+1$ (resp.\ $-1$) of the eigenvalue, 
$\abs{\lambda_j}<1$ (resp.\ $\abs{\lambda_{-j}}>1$). The point masses of $d\mu$ are 
precisely those Dirichlet points inside gaps with $\sigma_j =+1$. The set of allowed 
Dirichlet data are single points for closed gaps and a circle ($[y_j, z_{j+1}]\times 
\{-1,1\}$ glued at the ends) for open gaps. Thus, the totality is a torus of dimension 
$\ell =\#$ of open gaps. 

\begin{theorem}\lb{T10.4} If $\Delta$ has $\ell$ open gaps, then the subset of 
$\{\alpha_j\}_{j=0}^{p-1}\in\bbD^p$ which, when periodized, have discriminant 
$\Delta$ is a torus of dimension $\ell$. The map from these $\alpha$'s to the 
possible Dirichlet date is a bijection. 
\end{theorem} 

Critical to at least one understanding of this result is that the Carath\'eodory 
function, $F$, has a minimal degree meromorphic continuation to the genus $\ell-1$ 
hyperelliptic Riemann surface associated to $\sqrt{\Delta^2 -4}$. 

There is a natural symplectic form on $\bbD^p$ so that the real and imaginary 
parts of the coefficients of $\Delta$ are the set of integrals of a completely 
integrable system; this is described in Section~11.11 of \cite{OPUC2} and in 
\cite{Nen}. The associated flows include the defocusing AKNS flow.

%%%%%%%%%%%%%%%%%%%%%%%%%%%%%%%%%%%%%%%%%%%%%%%%%%%%%%%%%%%%%%%%%%%%%
\section{The Szeg\H{o} Mapping and the Geronimus Relations} \lb{s11}
%%%%%%%%%%%%%%%%%%%%%%%%%%%%%%%%%%%%%%%%%%%%%%%%%%%%%%%%%%%%%%%%%%%%%

Finally, we discuss a deep connection between OPRL and OPUC found by 
Szeg\H{o} \cite{Sz22a}. The map $z\mapsto z+z^{-1}$ maps $\bbD$ 
biholomorphically to $\bbC\cup\{\infty\}$ with a cut $[-2,2]$ removed. 
The map on the boundary, $e^{i\theta}\to 2\cos\theta$, is a two-to-one 
map of $\partial\bbD$ to $[-2,2]$ that induces a map from $\calM_{+,1} 
([-2,2])$ to those measures on $\partial\bbD$ which are invariant under 
complex conjugation. It is easy to see $\mu\in\calM_{+,1}(\partial\bbD)$ 
has such invariance if and only if its Verblunsky coefficients are real. 
Explicitly, $\rho$, a probability measure on $[-2,2]$, is associated to 
$\mu =\Sz(\rho)$, an even probability measure on $\partial\bbD$, via 
\begin{equation} \lb{11.1} 
\int f(x)\, d\rho(x) =\int f(2\cos\theta)\, d\mu(\theta)  
\end{equation}

Szeg\H{o} found the OPRL, $P_n$, for $\rho$ in terms of the OPUC, $\Phi_n$, 
for $\mu$: 
\begin{equation} \lb{11.2} 
P_n \biggl( z+\f{1}{z}\biggr) = [1-\alpha_{2n-1}(d\mu)]^{-1} z^{-n} 
[\Phi_{2n}(z) + \Phi_{2n}^*(z)] 
\end{equation} 
and used this to convert Szeg\H{o} asymptotics for OPUC (see \eqref{8.14}) 
to asymptotics for suitable OPRL. This asymptotics is often called 
{\it Jost asymptotics\/} in the discrete Schr\"odinger literature. 

Geronimus \cite{Ger46} found the relation between the Jacobi parameters  
$\{a_n,b_n\}_{n=1}^\infty$ for $\rho$ and the Verblunsky coefficients 
$\{\alpha_n\}_{n=0}^\infty$ for $\mu$ (with $\alpha_{-1}\equiv 1$): 
\begin{align} 
a_{n+1}^2 &= (1-\alpha_{2n-1})(1-\alpha_{2n}^2)(1+\alpha_{2n+1})  \lb{11.3} \\
b_{n+1} &= (1-\alpha_{2n-1}) \alpha_{2n}-(1+\alpha_{2n-1}) \alpha_{2n-2} \lb{11.4}  
\end{align} 

The map from $\alpha$ to $(a,b)$ is local, that is, changing a single $\alpha$ 
only changes a finite number of $a$'s and $b$'s. That is not true for the 
inverse. Scaled Chebyshev polynomials of the first kind has $a_1 =\sqrt2$, 
$a_n=1$ ($n\geq 2$), $b_n=0$, and the corresponding $\alpha_n\equiv 0$. 
Scaled Chebyshev polynomials of the second kind have $a_n\equiv 1$, 
$b_n =0$ (i.e., they differ at a single $a_n$), but have $\alpha_{2n}=0$ 
and $\alpha_{2n-1} = -1/(n+1)$. 

Still the inverse can be computed (\cite{Ger46}). Given $\{a_n, b_n\}_{n=1}^\infty$,  
define $\varphi_n^\pm$ by $\varphi_0 =0$, $\varphi_1 =1$, and for $n\geq 1$, 
\begin{equation} \lb{11.5} 
\varphi_{n+1}^\pm + a_{n-1}^2 \varphi_{n-1}^\pm + b_n \varphi_n^\pm 
= \pm 2 \varphi_n^\pm  
\end{equation}
By a Sturm oscillation theorem, $\{a_n,b_n\}_{n=1}^\infty$ are the Jacobi 
parameters of a measure supported on $[-2,2]$ if and only if $\varphi_n^+ 
>0$ and $(-1)^n \varphi_n^- >0$. The Verblunsky coefficients are given by 
\begin{alignat}{3} 
u_n &= \f{\varphi_{n+2}^+}{\varphi_{n+1}^+} \qquad 
& v_n &= -\f{\varphi_{n+2}^-}{\varphi_{n+1}^-} \lb{11.6}  \\
\alpha_{2n} &= \f{v_n -u_n}{v_n + u_n} \qquad  
& \alpha_{2n-1} &= 1-\tfrac12 \, (u_n + v_n) \lb{11.7}
\end{alignat} 

Recently, these mappings have been used by Denisov \cite{DenPAMS} and 
Damanik-Killip \cite{DamKil,Sim287} as a powerful tool in the study of 
discrete Schr\"odinger operators and of OPRL. For proofs and references, 
see Chapter~13 of \cite{OPUC2}.

\bigskip
%%%%%%%%%%%%%%%%%%%%%%%%%%%%%

\end{document}